%% file: main.tex
\documentclass{lmcs} 
\pdfoutput=1
\usepackage[utf8]{inputenc}

\usepackage{lastpage}
\lmcsdoi{19}{4}{16}
\lmcsheading{}{\pageref{LastPage}}{}{}%
{Dec.~23,~2022}{Nov.~28,~2023}{}

\usepackage{bm}
\usepackage{dsfont}

\usepackage{hyperref}
\usepackage{cleveref}
\usepackage{xspace}
\newcommand\mfbox{\fcolorbox{gray}{white}}
\usepackage{wrapfig}

\usepackage[]{algorithm2e}

\input{macros}


\begin{document}

\title[Density of LRS]{Computing the Density of the Positivity\texorpdfstring{\\}{ }Set for Linear Recurrence Sequences}

\author[E.~Kelmendi]{Edon Kelmendi\lmcsorcid{0000-0003-3100-1500}}
\address{Queen Mary University of London}
\email{e.kelmendi@qmul.ac.uk}

\begin{abstract}
  The set of indices that correspond to the positive entries of a sequence of numbers is called its positivity set. In this paper, we study the density of the positivity set of a given linear recurrence sequence, that is the question of how much more frequent are the positive entries compared to the non-positive ones. We show that one can compute this density to arbitrary precision, as well as decide whether it is equal to zero (or one). If the sequence is diagonalisable, we prove that its positivity set is finite if and only if its density is zero. Further, arithmetic properties of densities are treated, in particular we prove that it is decidable whether the density is a rational number, given that the recurrence sequence has at most one pair of dominant complex roots.

  Finally, we generalise all these results to symbolic orbits of  linear dynamical systems, thereby showing that one can decide various properties of such systems, up to a set of density zero. 
\end{abstract}


\maketitle

\section{Introduction}
\label{sec:intro}
Linear recurrence sequences (\lrs) are infinite sequences of rational numbers $\sq{u_n}$, whose every entry is a linear combination of the $k$ preceding entries. That is, a sequence that satisfies a recurrence relation:
\begin{align}
  \label{eq:rec rel}
  u_n = a_1u_{n-1}+\cdots+a_ku_{n-k},
\end{align}
for all $n>k$, where $a_1,\ldots, a_k$ are rationals and $a_k\ne 0$. The constants $a_1,\ldots, a_k$, and $u_1,\ldots,u_k$ uniquely identify the sequence. 

Firmly grounded as one of the fundamental families of finitely represented number sequences, they are ubiquitous in mathematics and computer science; their importance is evident. A basic object of study in modern number theory, they appear in the investigation of pseudo-random number generators, in cellular automata, as solutions of some Diophantine equations, as the number of $\mathbb F_q$-points on varieties, to name just a few examples. Furthermore, they are intrinsically related to linear dynamical systems, and the field of dynamical systems as a whole.

From another point of view, a linear recurrence sequence can be seen as a kind of restricted Turing machine, namely one that has a single loop inside which the variables are updated by a linear function. As such programs permeate any larger piece of software, verifying their correctness has become increasingly important in recent years. This motivation has driven further interest in algorithmic questions regarding these sequences. 

This field has been a rather active area of research in the past few decades --- a considerable body of work has amassed. The wide-scoped monograph~\cite{recseq} by Everest, van~der~Poorten, Shparlinski, and Ward is a place where one can find central results, their applications, as well as a taste of techniques that have proven useful. Here we recount only a brief summary of the theorems that are directly relevant to the present work.

We start with a basic question: What does the zero set of a linear recurrence sequence $\set{n\st u_n=0}$ look like? The wonderfully simple answer, provided in 1934 by Thoralf Skolem~\cite{skolem1934verfahren} using $p$-adic analysis, is that the zero set of a linear recurrence sequence is a finite union of arithmetic progressions and a finite set. In other words, the zero set is ultimately periodic. This theorem was soon after generalised to sequences of algebraic numbers by Mahler~\cite{mahler1935arithmetische}, and then later on by Lech, to sequences of members of any ring of characteristic zero~\cite{lech1953note}. An elementary proof of Skolem's theorem can be found in~\cite{hansel1985simple}, see also the discussion in Chapter 2.1 of~\cite{recseq}. Unfortunately, even though we know the form of zero sets, we do not know how to decide if it is empty. Every known proof of this result uses, in some way or other, $p$-adic analysis, resulting in a non-constructive argument. The question of whether one can decide if there exists some $n$, such that $u_n=0$, known as Skolem's problem, remains to this day, the central open problem for \lrs.

However, there are some partial results for sequences of low order\footnote{The order of the sequence is the smallest $k$ for which the sequence satisfies a recurrence like~\eqref{eq:rec rel}.}: With the help of Baker's theorem for linear forms in logarithms of algebraic numbers, Mignotte, Shorey, and Tijdeman \cite[Theorem 2]{84_distan_between_terms_algeb_recur_sequen}, and in parallel Vereshchagin \cite[Theorem 4]{vereshchagin}, proved that for sequences of order at most four, one can decide whether their zero set is empty. In the direction of hardness, Skolem's problem is known to be NP-hard~\cite{blondel02_presen_zero_integ_linear_recur}.

One can raise the same questions about the positivity set $\set{n\st u_n>0}$. This set, however, unlike the zero set, does not admit a clean description. In fact the positivity problem (is there some $n$ such that $u_n>0$) is more general than the Skolem problem. That is, there is a polynomial reduction from Skolem's problem to the positivity problem (with a quadratic increase in the order). The positivity problem is known to be decidable for \lrs of order at most five \cite{ouaknine13_posit_probl_low_order_linear_recur_sequen}, where Baker's theorem plays a crucial role again. In the direction of hardness, a decision procedure for the positivity problem for \lrs of order six would allow one to compute the homogeneous Diophantine approximation type of a large class of transcendental numbers \cite[Theorem 5.2]{ouaknine13_posit_probl_low_order_linear_recur_sequen}. Which suggests that such a procedure must come hand-in-hand with a deeper understanding --- than hitherto exists --- of Diophantine approximations of transcendental numbers.

Questions of asymptotic nature seem to be slightly more approachable. For example, one can decide if a sequence has infinitely many zeros \cite[Theorem 2]{berstel76_deux_des_suites}. The corresponding problem for the positivity set, \ie are there infinitely many $n$, for which $u_n>0$ is not known to be decidable, however. This problem is called the ultimate positivity problem\footnote{Ultimate positivity is the question: ``is it true that after some point every entry of the sequence is positive?''. If we ignore the zeros, ultimate positivity does not hold if and only if the negativity set is infinite, or the positivity set of $\sq{-u_n}$ is infinite.}.
In fact, as for positivity, a similar link to Diophantine approximations exists~\cite[Theorem 5.1]{ouaknine13_posit_probl_low_order_linear_recur_sequen}. Nevertheless, there is an important positive result: namely that the ultimate positivity problem is decidable for \emph{diagonalisable} \lrs~\cite{ouaknine14_ultim_posit_decid_simpl_linear_recur_sequen}. A sequence is diagonalisable if its characteristic polynomial, which for a sequence that satisfies~\eqref{eq:rec rel} is
\begin{align}
  \label{eq:char poly}
  x^k-a_1x^{k-1}-\cdots -a_{k-1}x-a_k, 
\end{align}
has no repeated roots. In fact, it is possible to go much further for diagonalisable sequences~\cite{almagor2021deciding}: One can decide any asymptotic $\omega$-regular property, even when the property itself is part of the input. For example, one can ask whether the sign pattern ``-+-'' occurs infinitely often in the sequence.

For the general case not much progress has been made however, it remains a long standing, difficult, open problem to decide anything about the positivity set of a general \lrs, in particular whether this set is empty, or whether it is finite. In the present paper, we prove that it is possible to decide some things about another notion of size of a subset of naturals: its \emph{density}. 

Recall that the density of a set $S\subseteq \nat$ is
\begin{align*}
  \lim_{n\to\infty}\frac {\lvert\set{1,2,\ldots, n}\cap S\rvert} n,
\end{align*}
where the vertical bars denote cardinality (note that the limit need not exist). The density is a notion used to measures how \emph{large} an infinite subset of natural numbers is. 

\begin{exa}
  Here is a trivial \lrs: $u_1=1$ and $u_n=-u_{n-1}$. Clearly its positivity set are the odd numbers, and its density is equal to $1/2$. 
\end{exa}
\begin{exa}
  \label{ex:2}
  It is possible to construct linear recurrence sequences that are equal\footnote{We have defined \lrs to be sequences of rational numbers, but one can define \lrs over larger rings, as is done for this example, where the ring is $\rel$. We chose the restriction to rationals for simplicity, although all the results of this paper can be proved for real algebraic numbers, at least.} to $\cos(n\theta)$, $n\in\nat$. If $\theta$ is a rational multiple of $\pi$, the positivity set of these sequences will have some rational density, if however $\theta$ is not a rational multiple of $\pi$ then, we will later see, that the density is equal to $1/2$.
\end{exa}
The density of the positivity set of any linear recurrence sequence always exists. This fact was proved by Bell and Gerhold \cite[Theorem 1]{bell05_posit_set_recur_sequen}, and is our principal starting point. With the exception of the paper above, to the best of our knowledge there is no other work that deals with the density of the positivity set. The paper \cite{berstel76_deux_des_suites} can however be interpreted as providing an algorithm to compute the density of the zero set. 

We now describe the results of this paper. The first one is of a qualitative nature:
\begin{thm}
  \label{th:den 1}
  There is a procedure that inputs a \lrs and decides whether the density of its positivity set is equal to 1. 
\end{thm}
The same procedure can be used to decide whether the density is equal to 0, after a trivial pre-processing step.

Bell and Gerhold have observed, by using an equidistribution theorem due to Weyl, a version of which can be found in Cassels's book \cite{cassels59_introd_to_dioph_approx}, that the density is equal to the Lebesgue measure of a certain set. We proceed along the same path and go further by constructing this set, for which it is necessary to explicitly describe the multiplicative relations among the roots of the polynomial in~\eqref{eq:char poly}. Afterwards, the problem is reduced to checking the emptiness of a semialgebraic set, which can be done using the decidability of the theory of real closed fields, \ie Tarski's algorithm. These tools have been successfully employed by Ouaknine, Worrell, and others, on a number of related problems, it is not surprising that they prove useful to bear on the problems of this paper as well.

We will show that this problem is both \np and co-\np hard, while the procedure in \autoref{th:den 1} runs in \pspace. When the order of the sequence is fixed, the complexity drops to \ptime. 

Although we do not yet know how to decide whether the sequence has infinitely many positive entries, we can decide whether there are \emph{many} of them, in the sense of having non-zero density, using \autoref{th:den 1}. Another point of view is that the question ``is the density 0?'' over-approximates the question ``is the positivity set finite?'', because a positive answer to the latter implies the same for the former. However, for the family of diagonalisable sequences, the implication becomes an equivalence --- the two questions are the same: 
\begin{thm} 
  \label{th:diagonalisable exact}
  In a diagonalisable sequence the positivity set is finite if and only if its density is zero. 
\end{thm}
\autoref{th:den 1} and \autoref{th:diagonalisable exact} together imply the main theorem of \cite{ouaknine14_ultim_posit_decid_simpl_linear_recur_sequen}, that ultimate positivity is decidable for diagonalisable \lrs. However the proof has the same ingredients, in particular a result on the growth of \lrs by Evertse, van~der Poorten and Schlickewei, which is based on a lower bound for sums of $S$-units, itself based on the deep ``subspace theorem'' of Schmidt.

The main theorem of this paper says that we can compute densities to arbitrary precision:
\begin{thm} 
  \label{th:compute}
  There is a procedure that inputs a \lrs $\sq {u_n}$ and a positive rational number $\epsilon\in\rat$, and computes some $\delta'$, such that $|\delta-\delta'|<\epsilon$, where $\delta$ is the density of the positivity set of $\sq{u_n}$. 
\end{thm}

The complexity is the same as for the density 1 problem; the problem is in \pspace in the description of the \lrs and $\lceil\epsilon^{-1}\rceil$, but it drops to polynomial time when the order of the sequence is fixed. 

The idea of the proof of \autoref{th:compute} is straightforward. We have to approximate the Lebesgue measure of a certain subset of the $d$-dimensional unit cube. To this end, we draw a grid of $N^d$ points and count how many of these points fall in the set. It then remains to prove that we can decide whether a given rational point is a member of the set, and to upper bound the error term. For the latter we use a result of Koiran~\cite{koiranil_approx}. We note that it is possible, instead of testing for every point whether it belongs to the set, to test it for fewer points that are picked randomly, resulting in a faster Monte-Carlo type algorithm.

Let us give a simple example that illustrates some of the ideas behind the theorems above.
\newpage
\begin{exa}
  \label{ex:3}
  Consider the following simple program:
  \begin{center}
    \mfbox{  
    \begin{minipage}{0.4\textwidth}
      \begin{algorithm}[H]
        x=0; y=6; z=4;

        \While{true}
        {
          $\begin{cases}
            x:=4x+3y\\
            y:=4y-3x\\
            z:=5z
          \end{cases}$

          \eIf{$y+z>0$}
          {
            {\color{red} Region A}
          }
          {
            {\color{blue} Region B}
          }
        }
      \end{algorithm}
    \end{minipage}}
\end{center}
where the assignments to the local variables $x,y,z$ are done in parallel.

It is not immediately evident from looking at this program that, for example, Region~A is entered infinitely often. The algorithm from \autoref{th:compute} can be used to conclude not only Region A is entered infinitely often, but that it is entered with frequency:
    \begin{align*}
      0.732279\ldots = \frac{\cos^{-1}(-2/3)} \pi. 
    \end{align*}
    At first sight, it might seem strange to see notions related to circles and triangles such as $\pi$ and $\cos$ appearing in the answer of a simple question about a simple program, but the reality is that only through them can we understand the program above. Let us explain the answer in more detail. The value of $y+z$ in the $n$-th iteration of the loop is clearly equal to
    
    \begin{align*}
      \begin{pmatrix}
        0 & 6 & 4
      \end{pmatrix}\cdot \begin{pmatrix}
        4 & -3 & 0\\
        3 & 4  & 0\\
        0 & 0 & 5
      \end{pmatrix}^n\cdot \begin{pmatrix}
        0 \\
        1 \\
        1
      \end{pmatrix}.
    \end{align*}
    Multiplying this quantity with $5^{-n}$ will not change its sign and we see that after the multiplication, the update matrix is a rotation in the first two coordinates:
    \begin{align*}
      \begin{pmatrix}
        4/5 & -3/5 & 0\\
        3/5 & 4/5 & 0\\
        0 & 0 & 1
      \end{pmatrix}^n=\begin{pmatrix}
        \cos n\phi & -\sin n\phi & 0\\
        \sin n\phi & \cos n\phi & 0\\
        0 & 0 & 1
      \end{pmatrix},
    \end{align*}
    where the angle $\phi$ is $\cos^{-1}(4/5)$. Multiplying this matrix with the row vector and the column vector above, we see that the variable $y+z$ in the $n$-th iteration of the loop has the same sign as $6\cos(n\phi)+4$. So now the question: with what frequency does the loop enter Region A? has been reduced to the question: for how many $n$ is $\cos(n\phi)>-2/3$? When the angle $\phi$ is not a rational multiple of $\pi$ (which is the case here), by Weyl's equidistribution theorem, $n\phi$ is uniformly recurrent modulo $\pi$, meaning that for any interval $I$ in $[0,\pi]$, the frequency with which $n\phi\mod \pi$ enters $I$ is proportional to the size of $I$ (that is length of the interval divided by $\pi$). As a consequence, since $\cos(n\phi)>-2/3$ if and only if $n\phi$ modulo $\pi$ belongs to the interval $[0,\cos^{-1}(-2/3)]$, the answer follows. 
\end{exa}

\autoref{ex:3} and \autoref{ex:2} show that density can be both a rational and an irrational quantity. Therefore, the algorithm in \autoref{th:compute} cannot \emph{a priori} be used to decide quantitative questions, such as whether the density is larger than some given rational. We give a partial result in this direction but leave the general case open:
\begin{thm} 
  \label{th:rational}
  There is a procedure that inputs a \lrs that has at most one pair of dominant complex roots, and decides whether the density of its positivity set is rational, and if it is, computes it exactly. 
\end{thm}

We also prove that when there are no (non-trivial) multiplicative relations among the dominant roots, the density is a \emph{period}, as defined by Kontsevich and Zagier~\cite{kontsevich01_period}. We note that conjectures by Kontsevich and Zagier, and of Grothendieck predict the transcendence degree of field extensions of $\rat$ generated by a finite set of intervals, but we do not pursue this conjectural direction further. 

Finally, we take a step back, and consider what makes the proofs of the theorems above work. One way of answering this question is to say that the \emph{sign sequence} of a \lrs is isomorphic (except in a set of density zero) to an $\omega$-word that belongs to a family of words, which we call toric words. Such words, we prove, have some pleasant properties, one of which is that we can compute the frequencies with which any given pattern appears. Taking this point of view, allows us to generalise 
Theorems \ref{th:den 1}, \ref{th:diagonalisable exact}, and \ref{th:compute},  to linear dynamical systems.

A \textbf{linear dynamical system} is given via a square $k\times k$ matrix $M$ with rational entries, and an initial point $x_0\in\rat^k$. Its orbit is the sequence of points
\begin{align}
  \label{eq:orbit}
  x_0, x_0 M, x_0 M^2, x_0 M^3, \ldots. 
\end{align}
In order to be able to illustrate an instance, suppose that $k=2$ and consider the following question.

\begin{wrapfigure}[19]{l}{.5\textwidth}
\begin{center}
\mfbox{  
  \begin{minipage}{0.4\textwidth}
    \includegraphics[width=\textwidth]{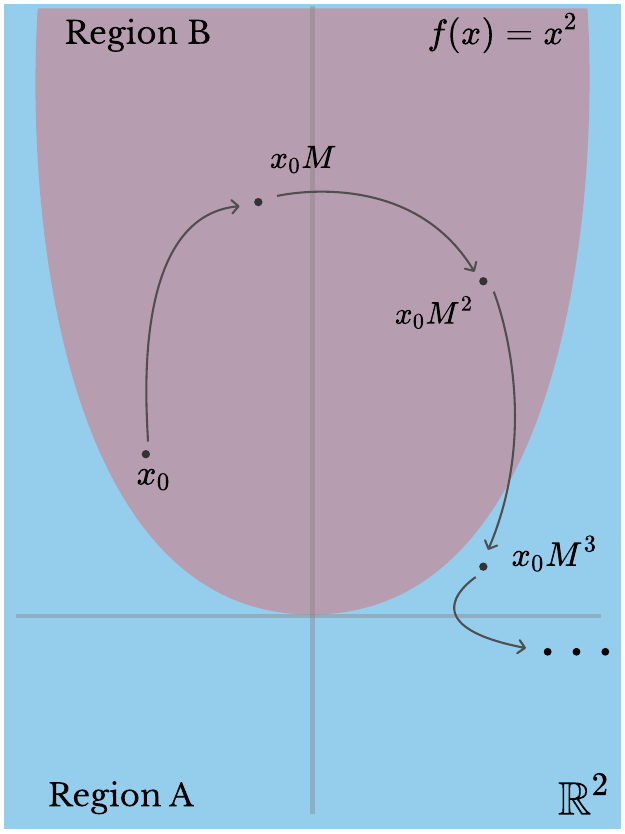}
  \end{minipage}
}
\end{center}
\end{wrapfigure}

\hfill\break
How frequently do members of the orbit in~\eqref{eq:orbit} enter Region B? More generally, suppose that we have a partition of the Euclidian space $\rel^k$ into semialgebraic sets $S_1,\ldots, S_\ell$. The latter are sets that one can define with polynomial inequalities, which we will define precisely later.  The orbit \eqref{eq:orbit} then defines an $\omega$-word $w$ over the alphabet $\set{1,2,\ldots,\ell}$, in the obvious way. This is a symbolic orbit of the dynamical system. One can then ask how frequently does a letter $b$ appear in $w$? In other words, what is the density of the subset of indices $n$ where $w_n=b$?  This and slightly more general questions can be answered by studying toric words; in the sense that there are analogues of Theorems \ref{th:den 1}, \ref{th:diagonalisable exact}, and \ref{th:compute}. The last section, \autoref{sec:toric}, is devoted to these questions.  

\vspace{1cm}

The rest of this paper is organised as follows. \autoref{sec:seqs} contains the principal definitions and generalities. \autoref{sec:non-degeneracy} is a technical section where we define a strong non-degeneracy condition and split the sequence into subsequences that satisfy it, as a pre-processing step for the algorithms that follow. \autoref{sec:den 1} deals with the density 1 problem, as well as the analysis for diagonalisable sequences. In the section that follows we give the procedure to compute the density. In the end, in \autoref{sec:rat}, we give the proof of \autoref{th:rational}, deciding when the density is a rational number.

A preliminary version of this paper appeared in \cite{licspaper}. 

\section*{Acknowledgment}
\noindent I am grateful to James Worrell for many helpful discussions.
\section{Sequences and Densities}
\label{sec:seqs}
A sequence $\sq{u_n}$ that satisfies a recurrence relation \eqref{eq:rec rel} for all $n>k$, but does not satisfy any linear recurrence with fewer terms, is called a \lrs of \textbf{order} $k$. The \textbf{characteristic polynomial} of such sequence is the polynomial \eqref{eq:char poly}, whose roots are, say
\begin{align*}
  \Lambda_1,\Lambda_2,\ldots,\Lambda_l,
\end{align*}
assumed to be distinct, with respective multiplicities $m_1\ldots,m_l$, where $1\le l\le k$. The sequence $\sq{u_n}$ can be written as a \textbf{generalised power sum} (see \cite[Section~1.1.6]{recseq}):
\begin{align}
  \label{eq:gps}
  u_n=\sum_{i=1}^lf_i(n)\ \Lambda_i^n,
\end{align}
where the polynomials $f_i$ have algebraic coefficients, $f_i\in\alg[x]$, and the degree of $f_i$ is  $m_i-1$. The converse also holds, any sequence $\sq{u_n}$ that can be written in the form \eqref{eq:gps} is a \lrs over algebraic numbers. The sequences whose roots all have multiplicity 1, \ie there are no repeated roots, are called~\textbf{diagonalisable} (or simple) sequences.

A \lrs is given by the numbers $a_1,\ldots, a_k$ and $u_1,\ldots,u_k$. From which, it is possible to compute descriptions of the constants in \eqref{eq:gps} in polynomial time in the bitlength of the input. By a \textbf{description} of an algebraic number we mean\footnote{There are other encodings of an algebraic number $\alpha$. Mostly one uses the fact that a number field $\rat(\alpha)$ is a vector space of finite dimension. For our purposes however, it is more convenient to define algebraic numbers by first-order formulas over the reals (defined below).} a first-order formula that defines it, typically this is the number's minimal polynomial together with intervals specifying where its real and imaginary parts lie. To compute the descriptions of the roots, one runs a root isolation algorithm on the characteristic polynomial (to compute the approximating intervals), see for example \cite{yap11} and \cite{06_real_roots}. Afterwards, for the computation of polynomials $f_i$, one solves a system of linear equations of polynomial size in the input. All this can be done in polynomial time. As a consequence, we assume that we have computed the descriptions of every constant in~\eqref{eq:gps}, and that the roots are ordered by their modulus, \ie
\begin{align*}
 |\Lambda_i|\ge |\Lambda_{i+1}|. 
\end{align*}
\textbf{Density} (also referred to as natural density, or asymptotic density in the literature) is a notion that measures how large a subset $S\subseteq\nat$ of natural numbers is. It is defined as:
\begin{align}
  \label{eq:density}
  \dens(S)\defeq\lim_{n\to\infty}\frac {\lvert\set{1,2,\ldots, n}\cap S\rvert} n, 
\end{align}
where by the vertical bars we denote the cardinality of the set. Not every set has a density; the limit might not exist. However they do have lower and upper density, which are defined by replacing limit with $\liminf$ and $\limsup$ respectively.
\begin{exa}
  Here is the density of some simple subsets of natural numbers. 
  \begin{enumerate}
  \item An (infinite) arithmetic progression, with common differences $d$, has density $1/d$. If the set $S\subseteq\nat$ is such that the difference between consecutive elements of $S$ is at most $d$, then the lower density of $S$ is larger than $1/d$. 
  \item The squares $\set{n^2\st n\in\nat}$ have density zero. To prove this, it suffices to observe that the cardinality of the squares in $\set{1,2,\ldots,n}$ is in $O(\sqrt n)$. 
  \item The primes have density zero due to the prime number theorem. 
  \end{enumerate}
\end{exa}

The principal object of study in this paper is the density of the \textbf{positivity set}:
\begin{align*}
  \dens\left(\set{n\st u_n>0}\right)
\end{align*}
of a given \lrs $\sq{u_n}$. Bell and Gerhold proved that it always exists:
\begin{thmC}[{\cite[Theorem 1]{bell05_posit_set_recur_sequen}}]
  The positivity set of any linear recurrence sequence has a density. 
\end{thmC}

The negativity set is just the positivity set of the sequence $\sq{-u_n}$ (which is plain, from \eqref{eq:gps} and the discussion above, that it can be computed). Therefore in the rest of this paper, we only deal with the density of the positivity set, which is simply referred to as the \textbf{density of the sequence}. 

We will make ample use of procedures for deciding the first-order \textbf{theory of real closed fields}, proved by Tarski \cite{tarski1951decision}. In this logic the atomic formulas are
\begin{align*}
  f(x_1,\ldots,x_n)\ge 0,
\end{align*}
where $f\in\intg[x_1,\ldots,x_n]$ is a polynomial with integer coefficients. The atomic formulas can be connected with Boolean connectives, and one is allowed to quantify over real numbers. Subsets of $\rel^n$ defined by such formulas are called \textbf{semialgebraic} sets. In the paper cited above, Tarski proved that there exists a procedure that inputs a first-order sentence and decides whether it is true when interpreted over the reals.

We can also interpret such formulas over the complex numbers instead of the reals, using the embedding of $\comp^n$ to $\rel^{2n}$, handling the real and imaginary parts individually. 

Note that our definition of descriptions of algebraic numbers is a simple formula in Tarski's logic. Other formulas that we will construct will be equally simple in the following sense: they will belong to the existential fragment, \ie formulas of the type
\begin{align*}
  \exists x_1\exists x_2 \cdots\exists x_n\qquad \Phi(x_1,\ldots, x_n),
\end{align*}
where $\Phi$ is quantifier-free. The complexity of this fragment is relatively low:
\begin{thm}[{\cite[Theorem 3.3]{Canny_1988} and \cite[Theorem 1.1]{renegar1992computational}, respectively}]
  \label{th:complexity of fo rel}
  The existential theory of reals is decidable in \pspace. When the number of variables is fixed, the complexity drops to \ptime\footnote{The \ptime upper bound holds for the full logic, when the number of variables is fixed, not only the existential fragment.}. 
\end{thm}
The theorems above expect the polynomials in the input to be written as a sequence of coefficients, each encoded in binary. Hence the exponents are assumed to be encoded in unary. 

\section{Strongly Non-Degenerate Subsequences}
\label{sec:non-degeneracy}
Let $P\in\nat$, and consider subsequences of the form: 
\begin{align}
  \label{eq:subseqs}
  \left\{\sq{u_{nP+\ell}}\st 0\le\ell <P\right\}. 
\end{align}
Each one is itself a \lrs (\cite[Theorem~1.3]{recseq}). One can easily observe this fact from the equality \eqref{eq:gps}: the roots of the subsequence are $\Lambda_i^P$ and the polynomials $f_i(nP+\ell)$ are multiplied by the constant $\Lambda_i^\ell$. 

The purpose of this section is a crucial preprocessing step that splits the sequence into subsequences \eqref{eq:subseqs}, for a particular $P\in\nat$, which we will compute. The subsequences have a number of properties (enumerated in a lemma at the end of the section) that make them more amenable. Effectively dividing the initial problem into easier sub-problems, we can recombine answers of the sub-problems to get the answer for the initial problem. For example, if we know the densities of the $P$ subsequences, then the density of the original sequence is equal to their sum divided by $P$. Or for the density 1 problem:  the original sequence has density 1 if and only if all the subsequences have density 1. 

In our case the period $P$ is a product:
\begin{align*}
  P\defeq P_1\cdot P_2,
\end{align*}
where $P_1$ comes from degeneracy, and $P_2$ from multiplicative relations among the roots. Let $N$ be the bitlength of the input and $k$ the order of the sequence, later in this section we will prove that $P$ will have the upper bound:
\begin{align}
  \label{eq:P ub}
  P\in 2^{\mathcal O(k^5\log\log N)}.
\end{align}
Before we give the definitions of the periods $P_1$ and $P_2$, let us first discuss the description of the roots $\Lambda_i^P$, as this is important for the complexity upper bounds when $P$ is large. Let $r\in\nat$, and let $z\in\alg$ be an algebraic number with description~$\phi(x)$ (\ie the formula $\phi(x)$ holds if and only if $x=z$). There are two ways to describe the number $z^r$:
\begin{enumerate}
\item The \emph{trivial way}: saying that there exists some $x$ such that $\phi(x)$ and
  \begin{align*}
    y=\underbrace{x\cdot x \cdots x}_{r\text{ times}}.
  \end{align*}
  Resulting in a constant increase on the number of variables, and a linear increase in $r$ on the size of the formula. 
\item The \emph{repeated squaring way}: saying that there exist a roughly $s:=\log r$ number of variables $x_1,\ldots,x_s$ such that
  \begin{align*}
    \phi(x_1)\text{ and } y=x_s\text{ and } x_{i+1}=x_{i}\cdot x_i,\, 1\le i\le s. 
  \end{align*}
  Resulting in a $\log r$ increase in both the number of variables and the size of the formula. 
\end{enumerate}
We will use both methods, depending on which complexity bound we want to derive. 
\begin{prop}
  \label{prop:formula size}
  For a given constant $P$, bounded by \eqref{eq:P ub}, the description of any $\Lambda_i^P$ can be computed in polynomial time. Furthermore, such a description grows both in the number of variables and in size by a term in $\mathcal O(k^5\log\log N)$.

  When the order of the sequence is fixed, the size of the description grows by a term in $\mathcal O(\log N)$ while the number of variables by a constant. 
\end{prop}
\begin{proof}
  Using the repeated squaring method results in a formula that grows both in size and in the number of variables by a $\log P$ term, hence the first statement of the proposition.

  When $k$, the order of the sequence is fixed however, it makes more sense to use the trivial way of constructing the formula, because this will result in a constant increase in the number of variables, and a linear in $P$ increase in the size of the formula. Since for fixed $k$, $P$ is in $\mathcal O(\log N)$ the second statement of the proposition follows. 
\end{proof}

Now we define $P_1$ and $P_2$, as well as show how to compute them. In the end of this section we summarise the properties that every subsequence $\sq{u_{nP+\ell}}$ has. 

\subsection{Period $P_1$}
\label{subsec:p1}

We begin with the standard notion of degeneracy. A \lrs is said to be \textbf{degenerate} if it has two distinct roots $\Lambda_i$ and $\Lambda_j$, whose ratio $\Lambda_i/\Lambda_j$ is a root of unity. One can test in \ptime whether a given sequence is degenerate by checking whether any of its ratios of distinct roots satisfies a cyclotomic polynomial of appropriate degree. Consult Section 3 in~\cite{yokoyama95_findin}. If the sequence is degenerate, taking the least common multiple of all the orders of roots of unity that can occur in this way, we get a quantity $P_1$, such that all the subsequences with period $P_1$ are either identically zero, or non-degenerate. The quantity $P_1$ is upper bounded only by a function in the order of the sequence: 
\begin{thmC}[{\cite[Theorem 1.2]{recseq}}]
  \label{th:subseqs}
  Let $\sq{u_n}$ be a \lrs of order $k$. Then there is a constant
  \begin{align*}
    M_k\in 2^{\mathcal O(k\sqrt{\log k})},
  \end{align*}
  such that for some $P_1\le M_k$, each subsequence
  \begin{align*}
    \sq{u_{nP_1+\ell}},
  \end{align*}
  $0\le\ell<P_1$ is either identically zero, or is non-degenerate. 
\end{thmC}

\subsection{Period $P_2$}
       
The definition of $P_2$ requires a little bit more work. We have assumed that the roots are ordered by their modulus: $|\Lambda_i|\ge |\Lambda_{i+1}|$, suppose that the first $j$ ones are dominant, \ie, 
\begin{align*}
  |\Lambda_1|=\cdots = |\Lambda_j| > |\Lambda_{j+1}|.
\end{align*}
Let $d$ be the maximal degree of the polynomials $f_1,\ldots,f_j$ from \eqref{eq:gps}, and suppose, without loss of generality, that it is exactly the polynomials $f_1,\ldots,f_m$ that are of degree $d$, for some $m\le j$. Define the normalised roots:
\begin{align*}
  \lambda_i\defeq \frac{\Lambda_i^{P_1}} {\left |\Lambda_i^{P_1}\right|}\qquad 1\le i\le m. 
\end{align*}
We are interested in the multiplicative relations: 
\begin{align*}
  \mathcal M(\lambda_1,\ldots,\lambda_m)\defeq\set{\vec r \in \intg^m\st \lambda_1^{r_1}\lambda_2^{r_2}\cdots \lambda_m^{r_m}=1}. 
\end{align*}
This set with addition forms a subgroup of $\intg^m$. Since the latter is a free abelian group with a basis of $m$ elements, by \cite[Theorem 7.3, Chapter I]{lang02_algeb} the subgroup~$\mathcal M$ is a free abelian group with some basis
\begin{align}
  \label{eq:basis}
  \vec b_1,\ldots,\vec b_v\in\intg^m,
\end{align}
where $v\le m$. Define
\begin{align}
  \label{def:p2}
  P_2\defeq 2\prod |b_{s,t}|,
\end{align}
where the product is taken over all $1\le s \le v$, and $1\le t\le m$, for which $b_{s,t}$ is nonzero. 
\begin{lem}
  \label{lem:p2 effective}
  The integer $P_2$ is effective. It can be computed in \pspace. When the order of the sequence is fixed, the computation can be performed in \ptime. 
\end{lem}
\begin{proof}
  We argue that we can compute the basis \eqref{eq:basis} and hence also $P_2$.

It follows from \cite[Theorem 1]{poorten77_multip_relat_number_field}, that there is an effective upper bound on the absolute value of the coordinates of the basis \eqref{eq:basis} of size:
\begin{align*}
  2^{\mathcal O(k^2)}\prod_{i=2}^m\log H(\lambda_i),
\end{align*}
where $H$ is the Mahler measure, defined as follows. For an algebraic number $z\in\alg$, with minimal polynomial
\begin{align*}
  a_0x^d+a_1x^{d-1}+\cdots + a_d=a_0(x-z_1)\cdots (x-z_d),
\end{align*}
we say that its Mahler measure is:
\begin{align*}
  H(z)\defeq|a_0|\prod_{i=1}^d\max\set{1,|z_i|} \le \sqrt{d}\max_{0\le i\le d}|a_i|,
\end{align*}
where the upper bound comes from~\cite[Lemma 1]{poorten77_multip_relat_number_field}. Using the fact that for any algebraic number $z\in\alg$ and $r\in\nat$, $H(z^r)=H(z)^r$, whose proof can be found in \cite[Chapter 3]{waldschmidt00_dioph_approx_linear_algeb_group}, via a straightforward computation, we can derive the following upper bound:
\begin{align}
  \label{eq:maxbst}
  \max_{1\le s \le v\atop{1\le t\le m}}|b_{s,t}| \in 2^{\mathcal O(k^3\log\log N)},
\end{align}
where $k$ is the order of the sequence and $N$ is the bitlength of the input. For any $\vec b\in\intg^m$ with the same upper bound, the assertion
\begin{align*}
  \vec b\in\mathcal M(\lambda_1,\ldots,\lambda_m),
\end{align*}
is an existential first-order formula of polynomial size in $N$, due to \autoref{prop:formula size}. Which means that by brute force, we can compute a basis \eqref{eq:basis} in \pspace by using the algorithm from~\autoref{th:complexity of fo rel}. When the order $k$ is fixed, the number of variables is constant. As a consequence of the second statement of~\autoref{th:complexity of fo rel}, in this scenario, the basis can be computed in \ptime.
\end{proof}

From the definition of $P_2$, \eqref{def:p2}, the estimate \eqref{eq:maxbst} and \autoref{th:subseqs}, one can derive the upper bound \eqref{eq:P ub}.

\subsection{Properties of the Subsequences}
Let $0\le\ell<P$, we list a number of properties of the subsequence
\begin{align}
  \label{eq: the subseq}
  \sq{u_{nP+\ell}},
\end{align}
which we assume is not identically zero. We start by replacing the dependent roots as follows.

The only case when the group $\mathcal M(\lambda_1,\ldots,\lambda_m)$ is trivial is when $m=1$, which implies that $\lambda_1=1$, because complex roots come as conjugate pairs (of the same multiplicity), and being a conjugate pair is a multiplicative relation (for algebraic numbers on the unit circle). In this case, every problem that we treat becomes trivial. Therefore suppose that $m>1$. Then there exists some member of the basis \eqref{eq:basis} --- say $\vec b_1$ without loss of generality --- that has at least two non-zero coordinates. By definition,
\begin{align*}
  \lambda_1^{b_{1,1}}\cdots \lambda_m^{b_{1,m}}=1.
\end{align*}
Suppose that $b_{1,m}\ne 0$. By using Euler's formula we see that we can write:
\begin{align}
  \label{eq:lambda m}
  \lambda_m=\varrho\lambda_1^{-b_{1,1}/b_{1,m}}\cdots \lambda_{m-1}^{-b_{1,m-1}/b_{1,m}},
\end{align}
where $\varrho$ is a $b_{1,m}$-th root of unity (and hence also a $P_2$-th root of unity, by definition of $P_2$), and at least one of the exponents $b_{1,1},\ldots,b_{1,m-1}$ is nonzero. Replacing $\lambda_m$ in the other equations, and continuing in this manner, making at most $v$ replacements, one for every member of the basis, we conclude that the set of indices $\set{1,\ldots,m}$ can be partitioned into the indices corresponding to the independent roots and depended roots of the form \eqref{eq:lambda m}, more precisely it can be partitioned into subsets:
\begin{itemize}
\item $I$ - a non-empty subset, with independent $\lambda_i$, \ie that do not have multiplicative relations among themselves, 
\item $D$ - a subset with dependent $\lambda_i$, \ie those that can be written in the form \eqref{eq:lambda m}, where in the right-hand side only members of $I$ appear, and there is a factor $\varrho$ which is a $P_2$-th root of unity (perhaps not primitive)\footnote{Here we see the reason behind the definition of $P_2$: In subsequences with the period $P_2$ we can directly write the dependent roots as a function of the independent ones; the factor $\varrho$ disappears because it is a $P_2$-th root of unity.}, and
\item $U$ - an empty set or a singleton containing some $i$ for which $\lambda_i=1$.
\end{itemize}
The reason why $U$ has cardinality at most $1$ is as follows. By the process described above, we cannot obtain more than one equation of the type $\lambda_i^r=1$, because among $\lambda_1,\ldots,\lambda_m$, the only root of unity that can appear is the number $1$. Indeed, if there were some complex $\lambda_i$ that is $r$-th root of unity, then its complex conjugate $\overline{\lambda_i}$ will also appear among the dominant roots $\lambda_1,\ldots,\lambda_m$ (with the same multiplicity), and $(\lambda_i/\overline{\lambda_i})^r=\lambda_i^{2r}=1$, meaning that the sequences $\sq{u_{nP_1+\ell}}$ are degenerate, a contradiction of \autoref{th:subseqs}.

Rearrange the the roots $\lambda_i$ such that for some $\eta$
\begin{align*}
  I=\set{1,\ldots,\eta},\qquad D=\set{\eta+1,\ldots,m-1},\qquad U=\set{m}. 
\end{align*}
The case when $D$ or $U$ is empty is omitted, as it can be treated in essentially the same way. It is convenient to define for all $i$, $1\le i\le m$:
\begin{align*}
  \alpha_i\defeq \lambda_i^{P_2}=\frac {\Lambda_i^{P}} {|\Lambda_i^{P}|},
\end{align*}
and the rationals $q_{i,j}\in\rat$, $i\in D$, $j\in I$, such that:
\begin{align*}
  \alpha_i=\prod_{j\in I}\alpha_j^{q_{i,j}}. 
\end{align*}
The generalised power sum form of the sequence \eqref{eq: the subseq} is:
\begin{align*}
  u_{nP+\ell}=\sum_{i=1}^{l}\Lambda_i^{\ell}f_i(nP+\ell)(\Lambda_i^P)^n.
\end{align*}
Dividing by $n^d|\Lambda_1^P|^n$ does not change the sign, where $d$ is the largest degree of polynomials multiplying the dominant roots. We get the sequence: 
\begin{equation}
  \label{eq:vn}
  \begin{aligned}
    v_n&\defeq\sum_{i=1}^mc_i\alpha_i^n+R(n)\\
    &=\sum_{i\in I}c_i\alpha_i^n+\sum_{i\in D}c_i\prod_{j\in I}\alpha_j^{q_{i,j}}+c_m+R(n),
  \end{aligned}
\end{equation}
where $c_i\in\alg$, and $R(n)$ is some residue that tends to zero polynomially, \ie
\begin{align}
  \label{eq:rn to zero}
  |R(n)|\in\mathcal O(n^{-\xi}),\text{ for some $\xi>0$}.
\end{align}
Furthermore there are no multiplicative relations among the roots $\alpha_i$, for $i\in I$, that is:
\begin{align}
  \label{eq:no multiplicative}
  \mathcal M(\alpha_1,\ldots,\alpha_{\eta})=\set{\vec{0}}.
\end{align}
A non-degenerate \lrs whose signs are the same as some sequence that can be written like $v_n$ above is what we call \textbf{strongly non-degenerate}. We summarise the properties of subsequences $\sq{u_{nP+\ell}}$.

\begin{lem}
  \label{lem:properties}
  For any $\ell$, $0\le\ell<P$, the following statements are true for the sequence $\sq{u_{nP+\ell}}$ that is not identically zero:
  \begin{enumerate}
  \item is non-degenerate,
  \item has finitely many zeros,
  \item its entries have the same sign as the entries of $\sq{v_n}$ defined in \eqref{eq:vn},
  \item the description of the algebraic numbers $c_i$, $\alpha_i$, and $q_{i,j}$ are of polynomial size, have polynomial many variables, and can be computed in \pspace,
  \item when the order of the sequence is fixed, the descriptions of the numbers above are of polynomial size, with a constant number of variables, and can be computed in \ptime. 
  \end{enumerate}
\end{lem}
\begin{proof}
  Property 1 comes from the fact that $P_1$ divides $P$ and \autoref{th:subseqs}. Any non-degenerate sequence that is not identically zero has finitely many zeros~\cite[Section 2.1]{recseq}, hence Property 2. The third property holds because we have obtained the sequence $\sq{v_n}$ by dividing with positive numbers.

  To see that Property 4 holds for the roots $\alpha_i$, first observe that $P$ can be computed in \pspace, as a consequence of \autoref{lem:p2 effective}, and the discussion in~\autoref{subsec:p1}. Then applying \autoref{prop:formula size} gives the wanted conclusion. One makes a similar argument  for the constants $c_i$. As for the rationals $q_{i,j}$, a combination of two facts is used. First, note that in the proof of \autoref{lem:p2 effective} the basis \eqref{eq:basis} is being computed in \pspace. Second, in the procedure that computes these rationals, described above, we do at most $v^2$ many replacements where $v\le m$ is the size of the basis. Each such replacement can be done in \ptime.

  For the last property, when $k$, the order of the sequence is fixed, the constant $P$ is in $\mathcal O(\log N)$, and it can be computed in \ptime, due to~\autoref{lem:p2 effective}. Note that in this case $P_1$ is constant. The property then follows by the same argument as for Property~4. 
\end{proof}

\autoref{ex:3} is not very interesting with respect to this section because after division by $5^n$ it is already in a strongly non-degenerate form. Here is a more suitable example.
\begin{exa}
  Let $\alpha$ be an algebraic number in the unit circle, for example:
  \begin{align*}
    \alpha\defeq \frac 3 5 + i\frac 4 5,
  \end{align*}
  and define:
  \begin{align*}
    \lambda_1\defeq\alpha^5,\qquad\lambda_2\defeq\alpha^3 \overbrace{\left(\frac {\sqrt{5}-1} 4 + i\frac {\sqrt{10+2\sqrt{5}}} 4\right)}^{\varrho}. 
  \end{align*}
  Then one can come up with a \lrs $\sq{u_n}$ over the real algebraic numbers\footnote{This \lrs is deliberately defined over the ring $\alg\cap \rel$ in order to keep the example small.}, whose characteristic polynomial, split into linear factors, is:
  \begin{align*}
    f(x)\defeq (x-\lambda_1)^2(x-\overline{\lambda}_1)^2(x-\lambda_2)^2(x-\overline{\lambda}_2)^2(x-1/2). 
  \end{align*}
  The sequence $\sq{u_n}$ in power sum form will look like:
  \begin{align*}
    u_n=(a_1+na_2)\ \lambda_1^n+(\overline a_1+n\overline a_2)\ \overline\lambda_1^n+(b_1+nb_2)\ \lambda_2^n+(\overline b_1+n\overline b_2)\ \overline\lambda_2^n+c\ 2^{-n},
  \end{align*}
  for some algebraic numbers $a_1,a_2,b_1,b_2,c$. First let us isolate the dominant terms, to this end, since $|\lambda_i|=1$, just divide the equality above by $n$, to get
  \begin{align*}
    a_2\lambda_1^n+\overline a_2\overline\lambda_1^n+b_2\lambda_2^n+\overline b_2\overline\lambda_2^n+R(n),
  \end{align*}
  where the remainder $R(n)$ tends to zero as $n\to\infty$. Since the ratio of $\lambda_i$ is a nonzero power of $\alpha$ (times $\varrho$ or $\varrho^{-1}$) it cannot be a root of unity. Hence the sequence is non-degenerate, \ie $P_1=1$. However, it is not strongly non-degenerate. Indeed, since $\varrho$ is a fifth root of unity, there is a multiplicative relationship between $\lambda_1$ and $\lambda_2$, which is
  \begin{align*}
    \lambda_1^3=\lambda_2^5. 
  \end{align*}
  So $(3,0,-5,0)$  and $(0,3,0,-5)$ form a basis of the subgroup of multiplicative relationships, hence $P_2$ in this case is equal to $450$. So we look at the strongly non-degenerate subsequences $u_{450n+\ell}$, where $\ell\in\set{0,\ldots, 449}$. Define $a_2':=a_2\lambda_1^\ell$, $b_2':=b_2\lambda_2^\ell$, and $\gamma_i=\lambda_i^{450}$. Then we have:
  \begin{align*}
u_{450n+\ell}=a_2'\gamma_1^n+\overline{a'}_2\overline\gamma_1^n+b_2'\gamma_2^n+\overline{b'}_2\overline\gamma_2^n+R(450n+\ell).
  \end{align*}
  Finally, since $\varrho$ is a primitive fifth root of unity, and therefore also a 450th root of unity we may write
  \begin{align*}
    \gamma_2=\gamma_1^{3/5},
  \end{align*}
  and with this replacement the equation above becomes:
  \begin{align*}
u_{450n+\ell}=a_2'\gamma_1^n+\overline{a'}_2\overline\gamma_1^n+b_2'\gamma_1^{3n/5}+\overline{b'}_2\overline\gamma_1^{3n/5}+R(450n+\ell).
  \end{align*}
  
\end{exa}

\section{The Density 1 Problem}
\label{sec:den 1}
In this section we prove that it is decidable whether the density of a given sequence is equal to $0$. The procedure expects a strongly non-degenerate sequence as input, \ie a sequence of the form in \eqref{eq:vn} with the properties that are listed in \autoref{lem:properties}. Suppose that we are given such a sequence and let $\delta$ be its density.

Note that the density of the \emph{negativity} set of the sequence (which is the same as the density of $\sq{-v_n}$) is equal to $1-\delta$, because the zeros $\sq{v_n}$ do not affect the density, being finitely many; a consequence of Property 2 in \autoref{lem:properties}. Hence the density of the sequence $\sq{v_n}$ is $0$ if and only if the density of $\sq{-v_n}$ is $1$. Thus the two problems, ``is the density 1?'' and ``is the density 0?'' are inter-reducible.

The argument for decidability of the density 0 problem is as follows. We define two open and measurable sets $\mathcal P$ and $\mathcal Q$ such that
\begin{align}
  \label{eq:nonempty}
  \mathcal P=\emptyset\qquad\Leftrightarrow\qquad\mathcal Q=\emptyset, 
\end{align}
and furthermore
\begin{align}
  \label{eq:mu}
  \mathcal Q\text{ is semialgebraic}\qquad\text{and}\qquad\delta=\mu(\mathcal P),
\end{align}
where $\mu$ denotes the Lebesgue measure. Being open sets, it follows that $\delta>0$ if and only if the semialgebraic set $\mathcal Q$ is nonempty, which can be decided, in particular because of \autoref{th:complexity of fo rel}. In this way decidability of the density 1 problem, \ie \autoref{th:den 1}, will follow from \eqref{eq:nonempty} and \eqref{eq:mu}, as well as the reduction from the density 1 to the density 0 problem.  

We proceed with the definitions of the sets $\mathcal P$ and $\mathcal Q$. Let $\torus$ be the unit circle, \ie the set of complex numbers $z\in\comp$, for which $|z|=1$.  Define the auxiliary functions $F$ and $G$ which are $v_n-R(n)$ but the roots $\alpha_i$ are replaced by variables; more precisely $F$ is a map from $[0,1]^\eta$ to the reals, and $G$ a map from $\torus^\eta$ to the reals, defined as: 
\begin{align*}
  &F(\vec\varphi)\defeq\sum_{i=1}^{\eta}c_i\exp(2\pi\ii\ \varphi_i)+\sum_{i=\eta+1}^{m-1}c_i\exp\left(2\pi\ii\sum_{j=1}^\eta q_{i,j}\varphi_j\right)+c_m,\\
  &G(\vec z)\defeq\sum_{i=1}^{\eta}c_iz_i\qquad\ \ \ \ \ \ \ \ +\sum_{i=\eta+1}^{m-1}c_i\prod_{j=1}^\eta z_i^{q_{i,j}}\qquad\qquad\ \ \ \ +c_m. 
\end{align*}
Now the sets $\mathcal P$ and $\mathcal Q$ are defined as:
\begin{align*}
  \mathcal P&\defeq\set{\vec\varphi\in [0,1]^\eta\st F(\vec\varphi)>0},\\
  \mathcal Q&\defeq\set{\vec z\in\torus^\eta\ \ \ \ \ \ \st G(\vec z)>0}.
\end{align*}
As one can obtain $\mathcal{P}$ by applying $\log z/2\pi\ii$ component-wise to elements of $\mathcal{Q}$, it is plain that $\mathcal{P}$ is non-empty if and only if $\mathcal{Q}$ is non-empty. Since $\mathcal{P}$ is open, it has non-zero measure if and only if it is non-empty. Furthermore, $\mathcal{Q}$ is semialgebraic, thus it only remains to show that $\delta=\mu(\mathcal P)$.

The proof follows closely the proof of the main theorem of~\cite{bell05_posit_set_recur_sequen}, and is crucially based on the following theorem, originally due to Weyl~\cite[Satz~4]{weyl16_gleic_von_zahlen_mod}, though we give a more modern reference from the book of Cassels. 
\begin{thmC}[{\cite[Theorem 1, page 64]{cassels59_introd_to_dioph_approx}}]
  \label{th:cassels}
  Let $\theta_1,\ldots,\theta_k,1\in\rel$  be linearly independent over $\rat$, and $S\subseteq [0,1]^k$ a measurable set, then
  \begin{align*}
    \dens\big(\set{n\st (n\theta_1\ \ \mathrm{mod}\ \ 1,\ldots,n\theta_k\ \ \mathrm{mod}\ \ 1)\in S}\big)=\mu(S).
  \end{align*}
\end{thmC}
It says that the fractional parts of $n\vec\theta$ fall in the set $S$ with frequency that is equal to the measure of the set $S$, in other words they are uniformly distributed in the $k$-dimensional cube. 

For $i\in\set{1,\ldots,\eta}$, define the arguments of the roots:
\begin{align*}
  \theta_i\defeq\frac{\log\alpha_i}{2\pi\ii}\in [0,1].
\end{align*}

Since there are no multiplicative relations among the $\alpha_1, \ldots, \alpha_\eta$, from \eqref{eq:no multiplicative}, we have that $\theta_1,\ldots,\theta_\eta,1$ are linearly independent over $\rat$. To see this, write $\alpha_i=\exp(2\pi\ii \theta_i)$ and observe that there are no multiplicative relations among the $\alpha_i$ if and only if there is no linear combination over $\rat$ of $\theta_i$ that is equal to an integer. As a consequence, we note that \autoref{th:cassels} is applicable to the tuple $\theta_1,\ldots ,\theta_\eta$.

The proof of $\delta=\mu(\mathcal P)$ is preceded by two lemmas. The first one says that the set of points that $F$ maps to $0$ has measure $0$.

\begin{lem}
  \label{lem:zero measure zero}
$\mu\big(\set{\vec\varphi\st F(\vec\varphi)=0}\big)=0.$
\end{lem}
\begin{proof}
  Since any generalised power sum is a \lrs over $\alg$ \cite[Section 1.1.6]{recseq}, the sequence
\begin{align*}
  \sq{F(n\vec\theta)}=\sq{v_n-R(n)}
\end{align*}
is a non-degenerate \lrs. As a corollary of the Skolem-Mahler-Lech theorem ~\cite[Section 2.1]{recseq}, this sequence has finitely many zeros, so
\begin{align*}
  \dens\big(\set{n\st F(n\vec\theta)=0}\big)=0.
\end{align*}
As noted above, we can apply \autoref{th:cassels} to $\vec\theta$, which implies
    \begin{align*}
    \dens\big(\set{n\st F(n\vec\theta)=0}\big)=\mu\big(\set{\vec\varphi\st F(\vec\varphi)=0}\big),
    \end{align*}
    where the set on the right-hand side is clearly measurable. Combining these two equations yields the statement of the lemma. 
  \end{proof}
  This lemma can also be proved without appealing to the Skolem-Mahler-Lech theorem, by directly showing that the set that is being measured has empty interior. 
  
The second lemma says that the indices in which the residue $R(n)$ is larger in absolute value than the dominating terms of the sequence, have upper density $0$. This means that it is only the dominant part that plays any role on the density $\delta$. Denote by $\udens$ the upper density (same as density except that the limit is replaced by $\limsup$): for all $S\subset \nat$,
\begin{align*}
  \udens(S)\defeq\limsup_{n\to\infty}\frac {\lvert\set{1,2,\ldots, n}\cap S\rvert} n.
\end{align*}
\begin{lem}
  \label{lem:R does not matter}
  $\udens\big(\set{n\st |F(n\vec\theta)|<|R(n)|}\big)=0$. 
\end{lem}
\begin{proof}
  For $\epsilon>0$, define:
  \begin{align*}
    \mathcal P_\epsilon&\defeq\set{\vec\varphi\in [0,1]^\eta\st |F(\vec\varphi)|\le\epsilon},\\
    \mathcal R_\epsilon&\defeq\set{n\in\nat\qquad \st |F(n\vec\theta)|\le \epsilon}.
  \end{align*}
  The residue $|R(n)|$ tends to zero as $n$ gets larger~\eqref{eq:rn to zero}, hence for all $\epsilon>0$,
  \begin{align}
    \label{eq:smaller than Reps}
    \udens\big(\set{n\st |F(n\vec\theta)|<|R(n)|}\big)\le\dens(\mathcal{R}_\epsilon). 
  \end{align}
  The set $\mathcal{R}_\epsilon$ has density as a consequence of \autoref{th:cassels}, also
  \begin{align*}
    \dens(\mathcal{R}_\epsilon)=\mu(\mathcal{P}_\epsilon)=\int_{[0,1]^\eta}\indic_{\mathcal{P}_\epsilon}d\mu,
  \end{align*}
  where by $\indic_{\mathcal{P}_\epsilon}$ we have denoted the indicator function of the set $\mathcal{P}_\epsilon$. Almost everywhere the function $\indic_{\mathcal{P}_\epsilon}$ tends to $\indic_{\mathcal{P}_0}$ as $\epsilon\to 0$, hence by Lebesgue's dominated convergence theorem~\cite[Theorem 16.4]{billingsleyprob} we have
  \begin{align*}
    \int_{[0,1]^\eta}\indic_{\mathcal{P}_\epsilon}d\mu\to\int_{[0,1]^\eta}\indic_{\mathcal{P}_0}d\mu=0,
  \end{align*}
  where the equality to zero comes from \autoref{lem:zero measure zero}. Since \eqref{eq:smaller than Reps} holds for all $\epsilon>0$, the statement of the lemma follows.
\end{proof}

One consequence of \autoref{lem:R does not matter} is that,
\begin{align*}
  \delta=\dens\big(\set{n\st v_n>0}\big)=\dens\big(\set{n\st F(n\vec\theta)>0}\big). 
\end{align*}
The density on the right-hand side is equal to $\mu(\mathcal{P})$ by again applying \autoref{th:cassels}.

Thus we have proved \autoref{th:den 1}, that it is possible to decide whether the density is equal to 0 (or to 1). The complexity of the procedure is in \pspace: the formula for non-emptiness of $\mathcal{Q}$ is of polynomial size due to Property 4 of \autoref{lem:properties}, and hence whether it is true can be decided in \pspace, \autoref{th:complexity of fo rel}.

The procedure runs in \ptime if the order of the sequence is fixed. This follows from Property 5 of \autoref{lem:properties} and \autoref{th:complexity of fo rel}.

Note that this lemma, \autoref{lem:R does not matter}, summarises the reason why we are able to decide certain properties of \lrs \emph{up to} a set of indices that has density zero. For a general \lrs it is rather difficult to understand for which indices $n\in\nat$, the dominant part $v_n$ is larger in absolute value than the absolute value of the remainder $|R(n)|$. Indeed this is the source of complications due to which we do not yet know whether the Skolem, positivity or ultimate positivity problems are decidable. It requires a deep understanding of certain arithmetic properties of the algebraic numbers $\alpha_i$. However, \autoref{lem:R does not matter} says that the indices $n$ for which the dominant part is smaller, form a subset of $\nat$ that has density zero, therefore in matters of density, these indices that are hard to understand have no effect.

\subsection{Complexity Lower Bounds}
It is possible to re-purpose the proofs of \cite{blondel02_presen_zero_integ_linear_recur} and \cite{ouaknine14_ultim_posit_decid_simpl_linear_recur_sequen} to show that the density 1 problem is both \np and co-\np hard. This indicates that the problem lies somewhere above these two classes, and is possibly \pspace-complete.

\begin{thm}
  The density 1 problem is \np-hard.
\end{thm}
\begin{proof}
  In essence, we will show that the proof of Blondel and Portier in \cite{blondel02_presen_zero_integ_linear_recur}, also implies the statement of the theorem. It works as follows.

  An instance of 3-\sat is a Boolean formula in variables $x_1,\ldots, x_n$ of the form:
  \begin{align}
    \label{eq:3sat}
    C_1\wedge C_2 \wedge \cdots \wedge C_m,
  \end{align}
  where each $C_i$ is the disjunction of exactly three terms, where a term is either $x_i$ or $\neg x_i$, $i\in \set{1,2,\ldots, n}$. The 3-\sat problem is \np-hard, and will be the problem we reduce from. The first reduction is into another problem, one about regular languages, which we describe now. 

  Fix a unary alphabet $\Sigma:=\set{a}$. A regular expression over this alphabet is built using the empty word $\epsilon$, words $a^n$, $n\in\nat$, union, and Kleene star. Here is an example of such a regular expression:
  \begin{align*}
    \epsilon\cup aaa \cup (\epsilon \cup aa)^*. 
  \end{align*}
  There is a polynomial reduction from 3-\sat to the problem that inputs such a regular expression and decides whether the language that it describes is different from the language~$a^*$. The reduction is as follows.

  Compute $p_1,\ldots, p_n$ the first $n$ prime numbers, which can be done in polynomial time (and they are all smaller than $n^2$). Define the function $h\st \nat\to\nat^n$ that maps
  \begin{align*}
    k\mapsto (k\ \mathrm{mod}\ p_1, k\ \mathrm{mod}\ p_2, \ldots, k\ \mathrm{mod}\ p_n),
  \end{align*}
  where by $k\ \mathrm{mod}\ p_i$ we denote the residue after dividing $k$ by $p_i$. Call a natural number $k$ a \textbf{code}, if and only if $h(k)\in\set{0,1}^n$.
  
  There is a regular expression $E_0$ such that the word $a^k$ belongs to the language $E_0$ describes, $L(E_0)$, if and only if $k$ is \emph{not} a code. That is the expression:
  \begin{align*}
    E_0 \defeq \bigcup_{i=1}^{n} \bigcup_{j=2}^{p_i-1} a^{j}(a^{p_i})^{*}. 
  \end{align*}
  Now let $C$ be one of the conjuncts in \eqref{eq:3sat}, and suppose that it involves the variables $x_r, x_s, x_t$. Consider a manner of setting bits $x_r, x_s, x_t$ such that the conjunct $C$ becomes false, \eg respectively $(x_r,x_s,x_t)=(0,1,0)$ makes $C=0$. Compute the smallest unique natural number $l$ such that
  \begin{align*}
    (l\ \mathrm{mod}\ p_r, l\ \mathrm{mod}\ p_s, l\ \mathrm{mod}\ p_t)=(0,1,0), 
  \end{align*}
  and the regular expression
  \begin{align}
    \label{eq:falsify}
    a^l(a^{p_r p_s p_t})^*. 
  \end{align}
  Let $E$ be the union of all such regular expressions (at most $8$ for each $C_i$) and of $E_0$. Denote by $L$ the language of $E$. Now by construction we have that the two following statements are equivalent for all $k\in\nat$:
  \begin{itemize}
  \item the word $a^k\not\in L$,
  \item $k$ is a code and the valuation $h(k)$ makes the formula \eqref{eq:3sat} true. 
  \end{itemize}
  Indeed, for the forward direction if $a^k$ does not belong to the language then it does not belong to $L(E_0)$ either, which means that it is a code, and it does not belong to the languages of expressions \eqref{eq:falsify} that encode valuations that falsify the conjuncts. The same argument can be used for the converse as well.

  By the Chinese reminder theorem we see that for any $v\in\set{0,1}^n$ there exists some $k$ such that $h(k)=v$. This then implies that the 3-\sat formula \eqref{eq:3sat} is satisfiable if and only if there is some $k\in\nat$ such that $a^k\not\in L$. Thus we have made the first reduction from 3-\sat.

  We observe one property of the language $L$ which we have just constructed. Define
  \begin{align*}
    p=\prod_{i=1}^np_i.
  \end{align*}
  By construction of $L$ we have that for all $k\in\nat$
  \begin{align}
    \label{eq:inf zeros}
    a^k\not\in L\qquad\Leftrightarrow\qquad a^{k+lp}\not\in L, \text{ for all }l\in\nat. 
  \end{align}
  Indeed, if $k$ is not a code then trivially $k+lp$ is not a code either, and if $k$ falsifies one of the conjuncts, then so does $k+lp$, by definition \eqref{eq:falsify}.

  Now we continue with the final reduction, from the problem about languages to the density 1 problem.

  From the regular expression $E$, construct in polynomial time a non-deterministic finite automaton $\mathcal A$ that recognises the language $L\setminus\set{\epsilon}$, and such that it has a unique initial and a unique final state. Suppose that its states are $\set{1,2,\ldots,t}$, where $1$ is the initial state and $t$ the final one. Let $M$ be the adjacency matrix of $\mathcal A$. Observe that the number $M^k_{i,j}$ is exactly the number of runs of length $k$ from state $i$ to state $j$. Then by construction, for all $k\in\nat$,
  \begin{align*}
    M^k_{1,t}\ne 0\qquad \Leftrightarrow\qquad a^k\in L. 
  \end{align*}
  The sequence $\sq{M^n_{1,t}}$ is in fact a \lrs whose every entry is non-negative. This \lrs has a zero if and only if the 3-\sat instance is satisfiable. From \eqref{eq:inf zeros}, if this \lrs has a zero then it has infinitely many of them, which fall on an infinite arithmetic progression with common differences at most $p$. In this case the density of the positivity set is $<1$, otherwise, if the sequence has no zeros, the density is equal to $1$. It follows that the density is not equal to one if and only if the 3-\sat instance \eqref{eq:3sat} is satisfiable. 
\end{proof}
\begin{thm}
  The density 1 problem is co-\np-hard.
\end{thm}
\begin{proof}[Proof Sketch]
  This lower bound follows immediately from \cite[Section 5]{ouaknine14_ultim_posit_decid_simpl_linear_recur_sequen}, so we give only a sketch.

  Consider the following problem. Given a polynomial $f\in\rat[x_1,\ldots, x_n]$ of degree at most $4$, decide whether there are real numbers $x_1,\ldots, x_n$ such that
  \begin{align*}
    f(x_1,\ldots, x_n)=0. 
  \end{align*}
  This problem, known as 4-\feas (for feasibility) is \np-hard, see for example \cite[Page~104, Theorem~1]{blum1998complexity}. The complement decision problem, \ie where one inputs a polynomial $f$ as above and one has to decide whether \emph{for all} real numbers $x_1,\ldots, x_n$ we have
  \begin{align}
    \label{eq:ge 0}
    f(x_1,\ldots, x_n)\ge 0,
  \end{align}
  is then co-\np-hard. By dividing the non-constant terms of $f$ with a certain integer that can be computed in polynomial time from $f$, we construct a different polynomial $f'$ such that \eqref{eq:ge 0} holds if and only if for all real $x_1,\ldots, x_n$ in the closed unit interval $[0,1]$, we have
  \begin{align}
    \label{eq:ge 0'}
    f'(x_1,\ldots, x_n)\ge 0. 
  \end{align}
  This problem is reduced in polynomial time to the ultimate positivity for \lrs in~\cite{ouaknine14_ultim_posit_decid_simpl_linear_recur_sequen}. The idea is to construct algebraic numbers $\lambda_1,\ldots, \lambda_n$ that lie on the unit circle, such that for all $i\in\set{1,\ldots,n}$ we have
  \begin{align*}
    \left\{\left(\lambda_i^k+\overline{\lambda}_i^k\right)^2\st k\in\nat\right\}\text{ is dense in }[0,1],
  \end{align*}
  and furthermore there are no multiplicative relations among the $\lambda_i$, and the expression in parenthesis is a linear recurrence sequence, with rational entries. Then we consider the sequence
  \begin{align*}
    f\left((\lambda_1^k+\overline\lambda_1^k)^2,\ldots, (\lambda_n^k+\overline\lambda_n^k)^2\right), k\in\nat,
  \end{align*}
  which is a \lrs; denote it by $\sq{u_n}$. Since the set of $x_1,\ldots, x_n$ in the $n$-cube $[0,1]^n$ for which $f(x_1,\ldots, x_n)<0$ is open (denote it by $X$), it follows that \eqref{eq:ge 0'} does not hold if and only if $u_n$ has infinitely many negative entries.

  It is possible, via the methods described in the beginning of the section, to conclude that by construction of $\sq{u_n}$, and \autoref{th:cassels}, the density of negative entries is equal to $\mu(X)$, the Lebesgue measure of $X$. Hence the density of the positive entries is equal to $1$ if and only if \eqref{eq:ge 0'} holds. The theorem follows. 
\end{proof}

Since in the beginning of this section we saw that an upper bound for the density 1 problem is \pspace, and the indications from the two theorems above are that a matching lower bound might exist, a search in this direction is interesting for the future. 

\subsection{The Case of Diagonalisable Sequences}
If the given \lrs has only finitely many positive entries then the density of the sequence is $0$. The converse, however, does not always hold, as it can be seen from the following example:
\begin{exa}
  One can construct an \lrs $\sq{w_n}$ that is equal to
  \begin{align*}
    w_n\defeq \frac n 2 \lambda^n + \frac n 2 \overline\lambda^n+(1-n),
  \end{align*}
  where $\lambda\in\torus$ is some algebraic number in the unit circle, that is not a root of unity. Let $\theta=\log\lambda/2\pi\ii$. Then, by writing $\lambda^n=\cos(2\pi n\theta)+\ii\sin(2\pi n\theta)$, we see that
  \begin{align*}
    w_n>0\qquad\Leftrightarrow\qquad \cos(2\pi n\theta)>1-\frac 1 n.
  \end{align*}
  The sequence $\sq{w_n}$ has infinitely many positive entries \cite[Proposition 4.1]{almagor2021deciding}. (This can be shown by appealing to Dirichlet's theorem \cite[Chapter~2, Theorem~1]{lang95_introd_dioph_approx}, and considering the Taylor's expansion of cosine.)

However the density of the positive entries is $0$. Indeed if it had density $\delta>0$, then we could have chosen some $n$ large enough such that the interval of $\varphi$ for which $\cos(2\pi\ii\ \varphi)>1-1/n$, is smaller than $\delta$, at which point, by applying \autoref{th:cassels} one can derive a contradiction. The latter theorem is applicable because $\lambda$ is not a root of unity, which means that $\theta$ is irrational, by definition. 
\end{exa}

The direction ``density 0'' implies ``positivity set is finite'', does however hold for an important class of \lrs, namely the diagonalisable sequences. These are sequences $\sq{t_n}$ whose characteristic polynomial has no repeated roots, as a consequence of which, its generalised power sum is of the following form:
\begin{align*}
  t_n\defeq\sum_{i=1}^ka_i\Lambda_i^n,
\end{align*}
where $a_i$ are some algebraic constants and $\Lambda_i$ are the roots. 

\begingroup
\def\thetheorem{\ref{th:diagonalisable exact}}
\begin{thm}
  In a diagonalisable sequence the positivity set is finite if and only if its density is zero. 
\end{thm}
\endgroup
\begin{proof}
  We prove the contrapositive, \ie we show that if $\sq{t_n}$ has infinitely many positive entries then it also has positive density. Assume that the roots are ordered by modulus, \ie $|\Lambda_i|\ge |\Lambda_{i+1}|$, and assume that the first $j$ roots have maximal modulus. Write
  \begin{align*}
    t_n=\underbrace{\sum_{i=1}^ja_i\Lambda_i^n}_{D(n)}+\underbrace{\sum_{i=j+1}^ka_i\Lambda_i^n}_{r(n)}.
  \end{align*}
  Suppose that $|\Lambda_1|>1$, indeed if it is not, we can always multiply the sequence with $\sq{K^n}$ for $K\in\nat$ large enough, without changing the sign. Without loss of generality, we can also assume that the sequence is non-degenerate. 

The proof hinges on a lower bound on the growth of \lrs that was proved Evertse, and in parallel by van der Poorten and Schlickewei, using the subspace theorem. See the discussion in \cite[Section 2.4]{recseq} as well as the appendix of~\cite{fuchs20_growt_linear_recur_funct_field}. Applying this theorem to our case, we have that for all $\epsilon>0$ there exists some threshold $n_0\in\nat$ such that:
  \begin{align*}
    |D(n)|\ge |\Lambda_1|^{(1-\epsilon)n}\text{ for all $n\ge n_0$}. 
  \end{align*}
  Since $|r(n)|$ can be upper bounded by some $c|\Lambda|^n$, with $c\in\rel$ a constant, and $|\Lambda|<|\Lambda_1|$, it follows that we can pick some $\epsilon>0$ for which we know that there exists some $n_0\in\nat$ such that:
  \begin{align*}
    |D(n)|>|r(n)|\text{ for all $n\ge n_0$}.
  \end{align*}
  This is a stronger version of \autoref{lem:R does not matter}, signifying that asymptotically the sign depends only on that of the dominant terms\footnote{This inequality holds for general \lrs. The difference is that for diagonalisable \lrs, the dominant part $D(n)$ is easier to analyse.}. As a consequence of the inequality above, since the sequence $\sq{t_n}$ has infinitely many positive terms, so does the sequence~$\sq{D(n)}$.

We sketch the rest of the proof. As in \autoref{sec:non-degeneracy} we can define the multiplicative relations among $\Lambda_1,\ldots,\Lambda_j$, and define a set $\mathcal{P}'$ analogous to the set $\mathcal{P}$, defined in the previous page. One can then prove that the set $\mathcal P'$ is open and furthermore it is non-empty as a consequence of the fact that $\sq{D(n)}$ has infinitely many positive entries. Non-emptiness implies that $\mathcal P'$ has non-zero measure, and finally, by applying \autoref{th:cassels}, one concludes that the density of the sequence is positive. 
\end{proof}
The algorithm that we have presented in this section is not the same, but it is quite similar to the algorithm of~\cite{ouaknine14_ultim_posit_decid_simpl_linear_recur_sequen} for deciding ultimate positivity for diagonalisable sequences. We have shown that this algorithm can be used for deciding a different problem, namely whether the density of the sequence is zero, and that when the sequence is diagonalisable, the density 0 question is equivalent to the question of whether the sequence has only finitely many positive entries. The complexity lower bound of \cite[Section 5]{ouaknine14_ultim_posit_decid_simpl_linear_recur_sequen} applies to our case as well.
\newpage
\section{Computing the Density}
\label{sec:comp}
One method of approximating the density $\delta$, which is the same as approximating the volume $\mu(\mathcal P)$ of the set $\mathcal P$ is conceptually simple: draw a grid and count the points that belong to $\mathcal P$. We summarise this in the picture below.

\begin{wrapfigure}[12]{l}{.4\textwidth}
\begin{center}
  \mfbox{
  \includegraphics[width=.25\textwidth]{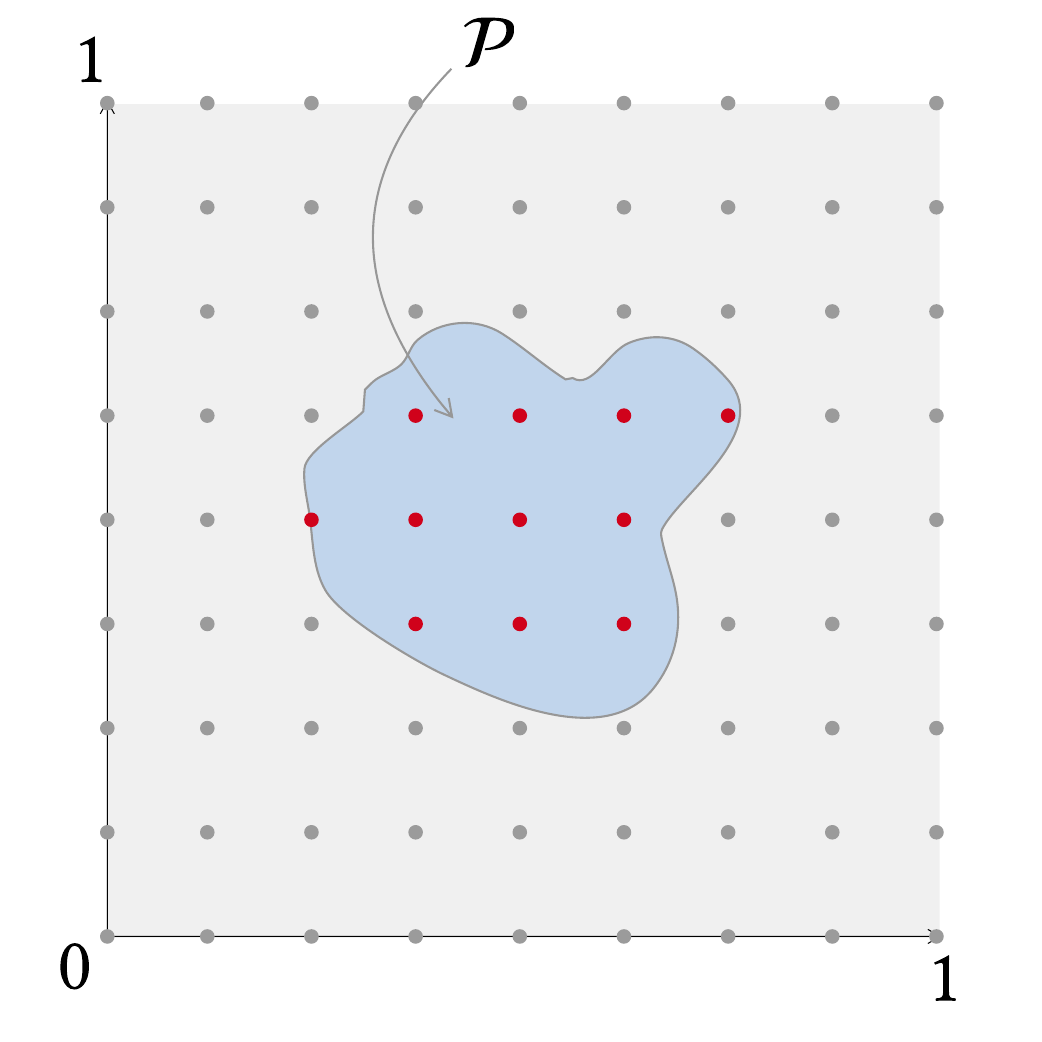}}
\end{center}
\end{wrapfigure}

From the grid of $M^\eta$ points (in the example $9^2$ points), we count how many are in $\mathcal P$, and denote this number by $C(M)$ (in the example this is equal to 11 red points). Since $\mathcal P$ is a measurable subset of the unit cube,
\begin{align*}
  \hskip 16em
  \frac{C(M)} {M^\eta} \to \mu(\mathcal P),
\end{align*}
as $M$ tends to infinity.

For this scheme to work, we need to be able to do two things. First, for any rational $\vec q\in [0,1]^\eta$, to be able to decide whether $\vec q\in\mathcal P$. And second, to be able to upper bound the quantity
\begin{align}
  \label{eq:error bound}
  \left|\frac {C(M)} {M^\eta} - \mu(\mathcal P)\right|,
\end{align}
by a function in $M$. We prove that both are feasible.

\begin{lem}
  \label{lem:rat points}
  Given any rational $\vec q\in [0,1]^\eta$, it is decidable whether $\vec q\in\mathcal P$.
\end{lem}
\begin{proof}
  Let $0\le k/n \le 1$ be a rational number. The complex number $\exp(2\pi\ii/n)$ is a primitive $n$-th root of unity, which we can isolate as a root of $x^n-1$. It follows that $\exp(2\pi\ii k/n)=\exp(2\pi\ii/n)^k$  is an algebraic number that we can easily define. Consequently the assertion $\vec q\in \mathcal P$, which is equivalent to $F(\vec q)>0$, is a first-order formula whose truth can be decided by Tarski's algorithm, \autoref{th:complexity of fo rel}. 
\end{proof}

For an upper bound on the error \eqref{eq:error bound}, we use the work of Koiran~\cite[Theorem 3]{koiranil_approx}. To introduce his theorem we need to define the parameter $\kappa(\mathcal P)$ first, and estimate it. 

Let $S\subseteq [0,1]^\eta$ be a measurable set, define $\kappa(S)$ to be the maximal number of connected components of the intersection $L\cap S$ where $L$ is an axis-parallel line. In other words, draw a line parallel to any one of the axes, and count how many times it goes in and out of the set. To estimate $\kappa(\mathcal P)$, in our case, this translates to fixing all but one parameter of function $F$ and counting how many times it will change its sign. More precisely, consider the function that we get by fixing all but one parameter of $F$, it will be of the form:
\begin{align*}
  H(\varphi)\defeq z_0\exp(2\pi\ii\ \varphi)+\sum_{i=1}^{\ell}z_i\exp(2\pi\ii\ r_i\varphi)+c,
\end{align*}
defined for $\varphi\in[0,1]$, where $\ell\le m$, $c$,$z_i$ are some algebraic numbers, and $r_i$ are taken among the $q_{i,j}$, $\eta<i<m$, $1\le j\le \eta$. The nature of the constants is such that $H$ is a real-valued function. How many times does $H$ change its sign in its domain $[0,1]$? By continuity, the answer is upper bounded by the number of zeros of $H$, which we will estimate. To this end, let $r_i=a_i/b_i$, for co-prime integers $a_i,b_i$, and define
\begin{align*}
  b\defeq \mathrm{lcm}\set{b_1,\ldots ,b_{\ell}}.
\end{align*}
Set $a'_i\in\nat$ to be such that $r_i=a_i'/b$.
\begin{lem}
  The function $H$ has at most
  \begin{align*}
    \hat q\defeq\max\set{b,a'_{1},\ldots,a'_{\ell}}
  \end{align*}
  zeros in the unit interval $[0,1]$. 
\end{lem}
\begin{proof}
  We can write $H$ as
  \begin{align*}
    z_0\left(\exp(2\pi\ii\ \varphi/b)\right)^b+\sum_{i=1}^{\ell}z_i\left(\exp(2\pi\ii\ \varphi/b)\right)^{a'_i}+c,
  \end{align*}
  which is a polynomial of degree at most $\hat q$, and hence can have at most that many zeros.
\end{proof}
Having estimated thus the parameter $\kappa(\mathcal P)$, we have the following upper bound on the error:
\begin{thmC}[{\cite[Theorem 3]{koiranil_approx}}]
  \label{thm:koiran}
  For all $M\in\nat$,
  \begin{align*}
    \left|\frac{C(M)} {M^\eta} - \mu(\mathcal P)\right|\le \frac {\eta\kappa(\mathcal P)} M\le\frac {\eta\hat q} M. 
  \end{align*}
\end{thmC}

Now \autoref{th:compute} follows from \autoref{lem:rat points} and \autoref{thm:koiran}. Indeed if we want to compute the density $\delta$ up to precision $\epsilon$, it suffices to choose $M\ge \eta\hat q/\epsilon$, then for every member of
\begin{align}
  \label{eq:grid}
  \left\{\left(\frac {k_1} M, \ldots, \frac {k_\eta} M\right)\st 0\le k_i\le M, 1\le i\le \eta\right\},
\end{align}
test whether it is in $\mathcal P$, and in this way compute the quantity $C(M)/M^\eta$ which by the proposition above is guaranteed to differ from the density by no more than $\epsilon$.

Even though $M$ is exponential in the input, by using the repeated squaring way of expressing the exponents in the formulas, as in~\autoref{sec:non-degeneracy}, it is possible to construct formulas of polynomial size for testing whether points of the grid \eqref{eq:grid} belong to $\mathcal P$. In particular to define $\exp(2\pi\ii/M)$, the formula says that it is a root of $x^M-1$ (which is of polynomial size), and that both the real and imaginary parts are positive and minimal. It follows that the algorithm for approximating the density is making exponentially many calls to a \pspace algorithm (due to \autoref{th:complexity of fo rel}), each of which is used to decide whether to increment a counter that is upper bounded by $M^\eta$. Hence this algorithm is running in \pspace on $\lceil\epsilon^{-1}\rceil$ and $N$, the bitlength of the description of the sequence. A similar analysis yields a \ptime upper bound in $N$ and $\lceil\epsilon^{-1}\rceil$ when the order of the sequence is fixed.

Instead of testing whether \emph{every} point in the grid belongs to $\mathcal P$, intuitively, we could test it for a smaller number $M'<M$, but choose the points uniformly at random.  This is the Monte-Carlo integration method \cite{koiranil_approx}. It results in a number of points in the set $C'(M')$ for which it is known that for all $\epsilon>0$,
\begin{align*}
  \frac 1 {M^\eta}\left|C'(M')-C(M)\right| \le \epsilon
\end{align*}
holds with probability at least $1-2e^{-2M'\epsilon^2}$. This can be demonstrated using Hoeffding's inequality.

\section{When is the Density a Rational Number?}
\label{sec:rat}
We have proved that it is possible to decide whether the density of a given sequence is equal to 0, or to 1. Can we also decide whether the density is larger than some $q\in\rat$? The approximating scheme of the previous section is \emph{a priori} of no help: since it might be the case that it outputs the estimates $\delta_1,\delta_2,\ldots$, for error bounds $\epsilon_1>\epsilon_2>\cdots$ such that $q$ belongs to all intervals $(\delta_i-\epsilon_i,\delta_i+\epsilon_i)$. A natural approach to tackling this difficulty is to ask whether the density itself is an irrational number. If the density is irrational then it has some $\epsilon$-neighbourhood which does not contain $q$, which means that for a sufficiently small $\epsilon_i$, $q$ does not belong in the interval $(\delta_i-\epsilon_i,\delta_i+\epsilon_i)$. 

In this section we report some progress of this direction. We begin by showing that when there are no non-trivial multiplicative relations among the roots, density is a \emph{period} as defined by Kontsevich and Zagier~\cite{kontsevich01_period}, \ie an integral of an algebraic function over a semialgebraic set. Afterwards, we prove that when there is at most one pair of dominant complex roots, it is decidable whether the density is rational, in which case we can compute it exactly.
\subsection{Density as a Period}
The complex roots of a sequence $u_n=\sum f_i(n)\Lambda_i^n$ come in conjugate pairs. Furthermore if $\Lambda_j=\overline{\Lambda_i}$ then also $f_j(n)=\overline{f_i(n)}$. See~\cite[Proposition 2.13]{halava2005skolem} for a proof. The multiplicative relations due to complex conjugacy, \ie $\lambda_j\lambda_i=1$, where $\lambda_i=\Lambda_i/|\Lambda_i|$ is the normalised root, we call \textbf{trivial relations}. Here we study sequences that do not have any non-trivial multiplicative relations among the roots. Under this restriction, the function $F$, defined in \autoref{sec:den 1}, has the following form:
\begin{align*}
  F(\vec\varphi)=\sum_{i=1}^{\eta}c_i\exp(2\pi\ii\ \varphi_i)+\sum_{i=1}^{\eta}\overline{c_i}\exp(-2\pi\ii\ \varphi_i)+c_m,
\end{align*}
since the only dependent roots are the complex conjugates of the independent ones. Using Euler's formula, and a trigonometric identity, we see that this function can also be written as:
\begin{align*}
  F(\vec\varphi)=c+\sum_{i=1}^\eta r_i\cos\left(2\pi(\varphi_i+\tau_i)\right),
\end{align*}
where $c=c_m\in\rel$, $r_i=|c_i|$, and $\tau_i$ is the argument of $c_i$.

We proceed by getting rid of the translation by $\tau_i$. Define:
\begin{align*}
  F'(\vec\varphi)\defeq F(\vec\varphi-\vec\tau).
\end{align*}
Recall that $\mathcal P$ was defined as the set of $\vec\varphi$ for which $F(\vec\varphi)>0$, and observe that
\begin{align*}
  \mathcal P'\defeq\set{\vec\varphi \st F'(\vec\varphi)>0}=\mathcal P + \vec\tau.
\end{align*}
Since $\mathcal P'$ is obtained from $\mathcal P$ by a translation, they have the same measure. Furthermore, as a consequence of symmetry of cosine we have:
\begin{align*}
  \mu(\mathcal P') = 2^{\eta}\mu\big(\underbrace{\mathcal P'\cap [0,1/2]^\eta}_{\hat{\mathcal P}}\big).
\end{align*}
So the density of the sequence can be derived from the volume of $\hat{\mathcal P}$. We write the latter as a certain integral. To this end, define the set $\mathcal L$ as,
\begin{align*}
  \mathcal L \defeq \left\{\vec x\in[-1,1]^\eta\st c+\sum_{i=1}^\eta r_ix_i>0\right\}. 
\end{align*}
Observe that the function $\cos^{-1}(\vec x)/2\pi$, denoted $g(\vec x)$, is a continuously differentiable bijection from $[-1,1]^\eta$ to $[0,1/2]^\eta$, and that furthermore:
\begin{align*}
  g(\mathcal L)=\hat{\mathcal P}. 
\end{align*}
Denote by $g'$ the Jacobian of $g$, then a variable change (see \cite[Theorem 3-13]{spivak18_calcul_manif}) leads to:
\begin{align*}
  \mu(\hat{\mathcal P})=\int_{g(\mathcal L)}d\vec\varphi=\int_{\mathcal L}|\det g'|d\vec x=\frac 1 {(2\pi)^\eta}\int_{\mathcal L}\prod_{i=1}^\eta\frac 1 {\sqrt{1-x_i^2}}d\vec x. 
\end{align*}
From here it follows that $\mu(\mathcal P)$ is rational if and only if
\begin{align}
  \label{eq:is it rat}
  \int_{\mathcal L}\prod_{i=1}^\eta\frac 1 {\sqrt{1-x_i^2}}d\vec x\ \ \ \in\ \ \ \rat\ \pi^\eta. 
\end{align}
The class of numbers that can be expressed as integrals of algebraic functions over semialgebraic sets are known as \textbf{periods}~\cite{kontsevich01_period}. They contain all algebraic numbers, as well as their logarithms, and some transcendental numbers like $\pi$; they are exceedingly commonplace however not well understood.

We do not know how to decide \eqref{eq:is it rat}, but we point out to some work that might prove to be helpful. One is Conjecture~1 in \cite{kontsevich01_period}, that says that if one period has two different representations as integrals, one can obtain one from the other through three simple operations: additivity, change of variables and Stokes's formula. It is not clear however, even if the conjecture were to be true, how one can calculate a sequence of such operations. A more direct conjecture is one made by Grothendieck that predicts the transcendence degree of field extension of $\rat$ that are generated by a finite set of periods. See \cite{ayoub2014periods} for definitions and a discussion about these two conjectures. More seems to be known about the special case of curves \cite{huber2018transcendence}, but in this case, for our purposes, we can give a more satisfactory answer by simpler means. 

\subsection{One Pair of Dominant Complex Roots}
When there is at most one pair of dominant complex roots, we have $\eta=1$ and the function $F$ can be written as:
\begin{align*}
  F(\varphi)=c+r\cos(2\pi(\varphi +\tau)). 
\end{align*}
Clearly when $|c|\ge |r|$ the density is either $1$ or $0$ depending on the sign of $c$, so assume that $|c|<|r|$. As we explained above, we can do away with the translation by $\tau$ when solely interested in density, and furthermore we can restrict $\varphi$ to the interval $[0,1/2]$.

Since the sequence is non-degenerate, the ratio $\lambda/\overline{\lambda}$ is not a root of unity, which implies that $\varphi$ is not a rational number. In this case, the equidistribution theorem, (\autoref{th:cassels}), is applicable. As a consequence of that theorem, to calculate the density, it suffices to calculate the length of the interval in $[0,1/2]$ which includes all $\varphi$ for which:
\begin{align*}
  \cos(2\pi\varphi)>\frac {-c} r. 
\end{align*}
Depending on the sign of $-c/r$, the length of this interval is
\begin{align*}
  \text{either}\ \ \ \ \frac{\cos^{-1}(-c/r)}{2\pi}\ \ \ \text{ or }\ \ \ 1-\frac{\cos^{-1}(-c/r)}{2\pi},
\end{align*}
in both cases it is rational if and only if $cos^{-1}(-c/r)$ is a rational multiple of $\pi$. In the remainder of this section we prove that we can decide whether the inverse cosine of a real algebraic number is a rational multiple of $\pi$. 
\begin{prop}
  \label{lem:arccos}
  Given a real algebraic number $\alpha\in[-1,1]$ of degree $d$, it is decidable whether
  \begin{align*}
    \cos^{-1}(\alpha)\in\rat\pi.
  \end{align*}
\end{prop}
\begin{proof}
Clearly $\cos^{-1}(\alpha)$ is a rational multiple of $\pi$ if and only if there is a rational $a/b\in\rat$, $b>0$, such that $a\cos^{-1}(\alpha)=b\pi$. Which, in turn, holds if~and~only~if there exists $a/b\in\rat$ (possibly different), $b>0$, such that:
\begin{align*}
  \cos\left(a\cos^{-1}(\alpha)\right)=(-1)^b.
\end{align*}
To proceed we need the definition of the \textbf{Chebyshev polynomials of the first kind} of order $n$. These are univariate polynomials $T_n$, for $n\in\nat$ that are characterised by the equation:
\begin{align*}
  T_n(\cos \theta)\defeq\cos(n\theta).
\end{align*}
One can also define them via a recurrence relation. We see that:
\begin{align*}
  T_a(\alpha)=T_a(\cos\cos^{-1}(\alpha))=\cos(a\cos^{-1}(\alpha))=(-1)^b.
\end{align*}
As a consequence $\cos^{-1}(\alpha)$ is a rational multiple of $\pi$ if and only if there is some $n\in\nat$, such that $\alpha$ is a root of
\begin{align*}
  T_n(x)-1\qquad\text{or}\qquad T_n(x)+1. 
\end{align*}
The roots of these polynomials are straightforward to describe: 
\begin{obs}
  \label{ob:1}
  Let $n\in\nat$. All the roots of $T_n(x)+1$ and of $T_n(x)-1$ come from the set
  \begin{align*}
    \big\{\pm\cos(k\pi/n)\st 0\le k \le n\big\}.
  \end{align*}
\end{obs}
The proof of this observation follows plainly from the fact that for all $x\in\rel$ such that $|x|\le 1$, we have
  \begin{align*}
    T_n(x)=T_n(\cos\cos^{-1}x)=\cos(n\cos^{-1}x)
  \end{align*}
and the fact that we can write $-\cos(k\pi/n)$ as $\cos(k\pi/n+\pi)$. 

From \autoref{ob:1} and the discussion preceding it we conclude that $\cos^{-1}(\alpha)$ is a rational multiple of $\pi$ if and only if it is equal to $\pm\cos(k\pi/n)$ for some $k,n\in\nat$, $k\le n$. The numbers $\pm\cos(k\pi/n)$ are algebraic, indeed they satisfy the Chebyshev polynomial of order $n$, furthermore if $\gcd(k,n)=1$ then $\cos(2k\pi/n)$ is an algebraic integer of degree $\Phi(n)/2$ \cite[Theorem~1]{lehmer1933note}, where $\Phi$ is the Euler's totient function.

Now, since $\alpha$ has degree $d$, we take some $N\in\nat$ such that $\Phi(N)\ge 2d$.  By testing (with the algorithms from \autoref{th:complexity of fo rel} say) whether $\alpha$ is a root of any $T_n(x)\pm 1$, for $n\le N$ we can decide whether $\cos^{-1}(\alpha)$ is a rational multiple~of~$\pi$.
\end{proof}
\newpage
\section{Toric Words and Linear Dynamical Systems}
\label{sec:toric}
As was described in the introduction, the decidability results of the preceding sections can be generalised to orbits of linear dynamical systems (\lds), where the positivity set is replaced by the set of indices corresponding to the members of the orbit that belong to a given semialgebraic set. In this section we explain how to achieve this generalisation.

Consider the orbit of a given \lds on the Euclidian plane.

\begin{center}
  \mfbox{
  \includegraphics[width=.9\textwidth]{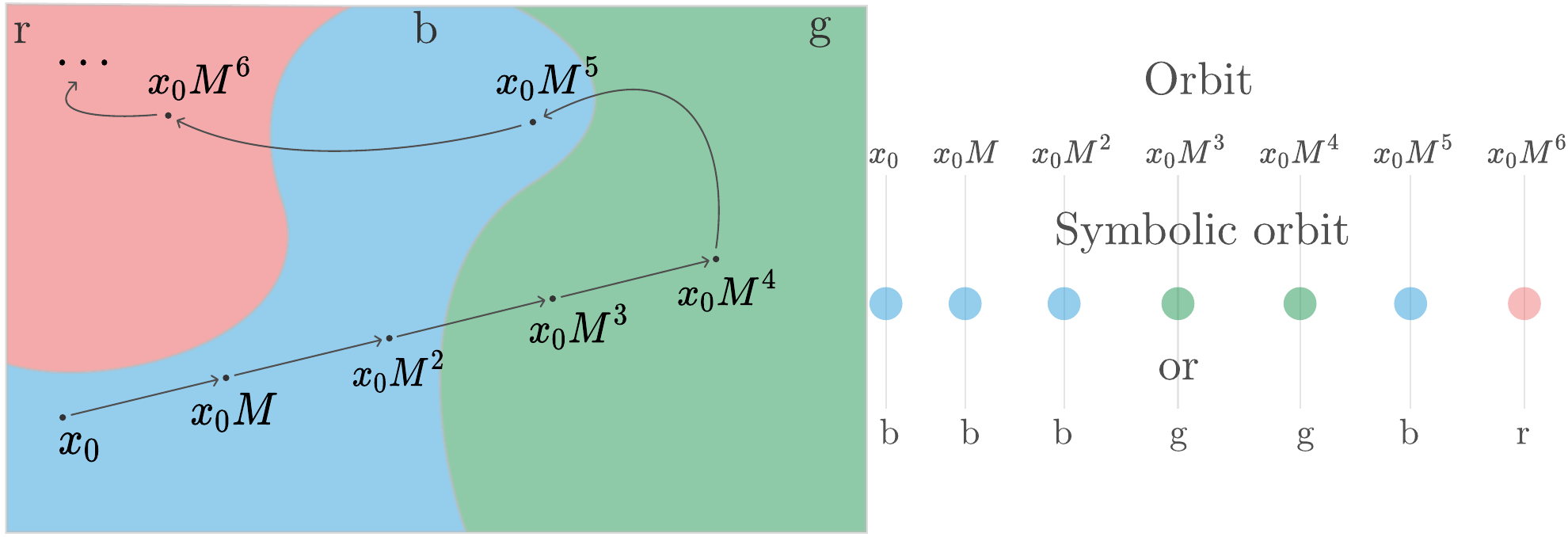}}
\end{center}

The latter is partitioned into some semialgebraic sets; in the example above in the red, blue, and green set. The orbit of the system is the sequence $x_0M^n$, $n\in\nat$, whereas the symbolic orbit (for the red, blue, and green subsets of $\rel^2$) retains only the information to which set the element of the orbit belongs. Hence it is an $\omega$-word over the alphabet $\Sigma:=\set{r,b,g}$. We want to compute how frequently some pattern occurs  in this word. More precisely, let $s=bbbggbr\cdots$ be the symbolic orbit, and $w\in\Sigma^*$ a finite word (or pattern) of length $|w|$, and say that $w$ \textbf{occurs} in $s$ in position $n$ if and only if
\begin{align*}
  s(n)s(n+1)\cdots s(n+|w|)= w,
\end{align*}
where by $s(n)$ we write the $n$th letter of $s$. We will show how to compute the density of the set:
\begin{align*}
  \set{n\st w\text{ occurs in $s$ in position $n$}},
\end{align*}
as well as decide whether it is equal to $0$ or $1$. We will call it the \textbf{density of the pattern}~$w$ in $s$, and denote it by
\begin{align*}
  \dens(w,s).
\end{align*}
This number gives rather precise (albeit asymptotic) information about the dynamics of the given \lds, namely it tells you in which set of the partition the system spends most its time.

How does this generalise the density of the positivity set of an \lrs? Let $k\in\nat$ and suppose that the \lrs $\sq{u_n}$ is given with the recurrence relation:
\begin{align*}
  u_n=a_1u_{n-1}+\cdots+a_ku_{n-k},
\end{align*}
where $a_k\ne 0$, and the first $k$ entries: $u_1,\ldots, u_k$. Denote by $M$ its companion matrix:
\begin{align*}
  M\defeq
  \begin{pmatrix}
    0 & \cdots & 0 & a_k\\
    \ & \      & \ & a_{k-1}\\
    \ & I      & \ & \vdots\\
    \ & \      & \ & a_1
  \end{pmatrix},
\end{align*}
where the block marked by $I$ is the $(k-1)\times (k-1)$ identity matrix. Denote by $\mathbf v$ the row vector $(u_1,\ldots,u_k)$. Then clearly we have for all $n\in\nat$,
\begin{align*}
  \text{$k$th coordinate of }\mathbf v M^n\text{ is equal to }u_{n+k}. 
\end{align*}
The partition of $\rel^k$ that we take is then the semialgebraic set
\begin{align*}
  S\defeq\set{(x_1,\ldots, x_k)\in \rel^k\st x_k>0},
\end{align*}
and its complement $\tilde S$. The symbolic orbit will be an $\omega$-word over a binary alphabet, where one letter (call it $p$) would imply that the corresponding entry is positive, while the other letter would imply that it is $\le 0$. Then the density of the positivity set of the \lrs $\sq{u_n}$ is just the density of the pattern $p$ in the symbolic orbit of $M$ for $(S,\tilde S)$. 

We give now the precise definitions. Let $d\in\nat$ and $\bm \lambda\in\torus^d$. The \textbf{orbit} of $\bm \lambda$ is the sequence: 
\begin{align*}
  \bm \lambda (n)\defeq (\lambda_1^n,\ldots,\lambda^n_d),\ \ n\in\nat. 
\end{align*}
Let $k\in\nat$, and  $S_1,\ldots, S_k\subset \torus^d$. The \textbf{symbolic orbit} of $\bm \lambda$ for $S_1,\ldots, S_k$ is an infinite word $s$ over the alphabet
\begin{align*}
  2^{\set{1,2,\ldots k}},
\end{align*}
defined as follows. For all $i\in\set{1,\ldots, k}$ and $n\in\nat$,
\begin{align*}
  i\in s(n)\qquad \Leftrightarrow\qquad \bm \lambda(n)\in S_i. 
\end{align*}
Same as in the example above (except that we have not assumed that the semialgebraic sets form a partition), the symbolic orbit is an abstraction of the orbit in which the only information we want to retain for a point $\bm \lambda(n)$ is to which sets $S_1,\ldots, S_k$ it belongs. 

Let $w_1\in\Sigma_1^\omega$ and $w_2\in\Sigma_2^\omega$. We say that $w_2$ is a \textbf{coarsening} of $w_1$ if there exists a map
\begin{align*}
 r\st \Sigma_1\to\Sigma_2,
\end{align*}
such that $w_2=r(w_1)$, applied letter-wise. In this case we also say that $w_1$ \textbf{refines} $w_2$. If both $w_1$ refines $w_2$ and vice-versa, we say that $w_1$ and $w_2$ are \textbf{isomorphic}; equivalently there exists an injective such map $r$, such that $w_2=r(w_1)$.

A \textbf{toric word} is then any word isomorphic to the symbolic orbit of some $\bm \lambda\in\torus^d$ for some $S_1,\ldots, S_k\subset \torus^d$, where the components of $\bm \lambda$ are algebraic numbers, and the sets $S_1,\ldots, S_k$ are semialgebraic.

Subfamilies of toric words have been studied in other contexts: For example, the special case $d=1$, $\lambda\in\torus$ is not a root of unity, and there is one set $S_1\subset \torus$ that is an interval, has a symbolic orbit that is a Sturmian word. Such symbolic orbits have been studied going back to Johann Bernoulli III (1744-1807), see the notes on Chapter~9 of~\cite{allouche2003automatic}. Muchnik et al. consider the case where $S_1,\ldots, S_k$ are open and disjoint and prove that in that case the symbolic orbits are almost periodic \cite[Section 4.3]{muchnik03_almos_period_sequen}. These words also have connections to extensions of \mso logic over $(\nat,<)$; however for the purposes of this paper we will be content with only showing a couple of closure properties of these words. 

It is convenient to assume that the semialgebraic sets $S_1,\ldots, S_k$ partition $\torus^d$, we can do this without loss of generality:
\begin{lem}
  \label{lem:partition}
  Let $\bm\lambda$ and $S_1,\ldots S_k$ be as above. There exists a partition of $\torus^d$ into semialgebraic sets $R_1,\ldots, R_h$ such that the symbolic orbit of $\bm\lambda$ for $S_1,\ldots, S_k$ is isomorphic to that of $\bm\lambda$ for $R_1,\ldots, R_h$.
\end{lem}
\begin{proof}
  Denote by $S_{\emptyset}$ the relative complement in $\torus^d$, of the union of $S_1,\ldots, S_k$. For all non-empty $J\subset\set{1,\ldots, k}$ define:
\begin{align*}
  S_J\defeq\bigcap_{i\in J}S_i-\bigcup_{i\not\in J}S_i. 
\end{align*}
Then the sets $S_J$ for $J\subset\set{1,\ldots, k}$ partition $\torus^d$, and they are furthermore semialgebraic. Enumerate these sets as $R_1,\ldots,R_h$. The symbolic orbit of $\bm \lambda$ for $R_1,\ldots, R_h$ is an infinite word over the alphabet $2^{\set{1,\ldots, h}}$. However, since the sets $R_1,\ldots, R_h$ partition $\torus^d$, the only letters that will appear in the symbolic orbit are the singletons $\set{i}$, $i\in\set{1, \ldots, h}$. 

It is not difficult to see now that the symbolic orbit of $\bm\lambda$ for $S_1,\ldots, S_k$ and that of $\bm\lambda$ for $R_1,\ldots, R_h$ are isomorphic, where the isomorphism depends on the particular enumeration that we have chosen. 
\end{proof}

We continue with the closure properties. Let
\begin{align*}
  w_1\in\Sigma_1^\omega,\qquad w_2\in\Sigma_2^\omega,
\end{align*}
be two infinite words over the respective alphabets $\Sigma_1,\Sigma_2$. The \textbf{product} of $w_1$ and $w_2$ is the word $w_1\times w_2$ over the alphabet $\Sigma_1\times \Sigma_2$ defined as
\begin{align*}
  (w_1\times w_2)(n)\defeq (w_1(n), w_2(n)). 
\end{align*}
\begin{prop}
  \label{prop:closed products}
  Toric words are closed under taking products. 
\end{prop}
\begin{proof}
    Let $w_1,w_2$ be two toric words, where $w_1$ is obtained from the symbolic orbit of $\bm\lambda\in\torus^d$ for semialgebraic sets $S_1,\ldots S_k$, and $w_2$ from that of $\bm\gamma\in\torus^e$ for $T_1,\ldots, T_e$. From \autoref{lem:partition}, we may assume that $S_1,\ldots, S_k$ and $T_1,\ldots, T_e$ partition $\torus^d$, respectively $\torus^e$. Then $\set{S_1,\ldots, S_k}\times \set{T_1,\ldots, T_e}$ partitions $\torus^{d+e}$ and furthermore those sets are semialgebraic. Now it is plain that the symbolic orbit of
\begin{align*}
  (\lambda_1,\ldots, \lambda_d,\gamma_1,\ldots,\gamma_e),
\end{align*}
for
\begin{align*}
  \set{S_1,\ldots, S_k}\times \set{T_1,\ldots, T_e}
\end{align*}
is isomorphic to the product of $w_1$ and $w_2$. 
\end{proof}
\begin{prop}
  \label{prop:closed under coarsenings}
  Toric words are closed under coarsenings. 
\end{prop}
\begin{proof}
  Let $w$ be isomorphic to the symbolic orbit of $\bm\lambda\in\torus^d$ for $S_1,\ldots, S_k$. We can assume that $S_1,\ldots, S_k$ partitions $\torus^d$ due to \autoref{lem:partition}. Then any coarsening of $w$ is isomorphic to the symbolic orbit of $\bm\lambda$ for $T_1,\ldots,T_l$, for some $t\le k$, where $T_i$ are made of unions of sets $S_i$, and are therefore semialgebraic. 
\end{proof}

A pleasant property of toric words, among others, is that we can decide whether the density of any given pattern that occurs in it is 0, or 1, as well as compute it to arbitrary additive precision. 

\begin{thm}
  There is a procedure for the following problem. Given as input:
  \begin{itemize}
  \item semialgebraic sets $S_1,\ldots, S_k\subset \torus^d$,
  \item algebraic numbers $\bm\lambda\in\torus^d$,
    
  \item a pattern $w\in 2^{\set{1,\ldots,k}^*}$,
  \end{itemize}
  decide whether the density of $w$ in the symbolic orbit of $\bm\lambda$ for $S_1,\ldots, S_k$ is zero. 
\end{thm}
\begin{proof}
  Due to \autoref{lem:partition}, we can assume that $S_1,\ldots, S_k$ partition $\torus^d$, and therefore that the pattern $w$ is a finite word over the alphabet $\Sigma:=\set{1,\ldots,k}$.

  Given $\bm\alpha,\bm\beta\in\torus^d$, we write $\bm\alpha\bm\beta$ for the vector:
  \begin{align*}
    (\alpha_1\beta_1,\ldots,\alpha_d\beta_d). 
  \end{align*}

  Denote by $m$ the length of the word $w$, and define the following semialgebraic set:
  \begin{align*}
    T_0\defeq \set{\bm\alpha\in\torus^d\st \bm\alpha\in S_{w(1)},\  \bm\alpha\bm\lambda(1)\in S_{w(2)},\ \ldots,\ \bm\alpha\bm\lambda(m-1)\in S_{w(m)}}. 
  \end{align*}
  The set $T_0$ characterises all points, starting from which, the orbit of $\bm\lambda$ moves among $S_1,\ldots, S_k$ in the next $m$ steps in pattern $w$. To rephrase this more precisely, denote by $s\in\set{1,\ldots,k}^\omega$ the symbolic orbit of $\bm\lambda$ for $S_1,\ldots, S_k$. And by $t$ the symbolic orbit of $\bm\lambda$ for $T_0,\tilde T_0$, where the latter is the relative complement of $T_0$ in $\torus^d$; then the infinite word $t$ is over the alphabet $\set{0,\tilde 0}$. By construction, the following statements are equivalent for all $n\in\nat$:
  \begin{enumerate}
  \item $w$ occurs in $s$ in position $n$,
  \item $t(n)=0$,
  \item $\bm\lambda(n)\in T_0$. 
  \end{enumerate}
  This equivalence implies that
  \begin{align}
    \label{eq:dens w s}
    \dens(w,s)=\dens(0,t),
  \end{align}
  where by $\dens(w,s)$ we have denoted the density of the pattern $w$ in $s$.
  
  The procedure computes the period $P$ as in \autoref{sec:non-degeneracy}, for the algebraic numbers $\bm\lambda$. This is a product of $P_1$ (the least common multiple of orders of roots of unity that one can obtain by taking ratios $\lambda_i/\lambda_j$, $i\ne j$), and $P_2$ that only depends on the multiplicative relations among the coordinates of $\bm\lambda$. We split the orbit of $\bm\lambda$ and the symbolic orbit $t$ into subsequences (subwords) by taking indices $nP+\ell$, $n\in\nat$, $0\le\ell<P$. Denote by $t_\ell$ such a subword, \ie
  \begin{align*}
    t_\ell(n)\defeq t(nP+\ell). 
  \end{align*}
  Clearly the density of $0$ in $t$ is positive (and hence also $\dens(w,s)>0$, due to (\ref{eq:dens w s})) if and only if there exists some $\ell$, $0\le\ell<P$, such that the density of $0$ in $t_\ell$ is positive. As the procedure tries to find such an $\ell$, we only need to show how to decide whether
  \begin{align*}
    \dens(0,t_\ell)>0,
  \end{align*}
  for some fixed $\ell$, $0\le\ell<P$.

  Define
  \begin{align*}
    T_{\ell}\defeq\set{\bm\alpha\in\torus^d\st \bm\lambda^{\ell}\bm\alpha\in T_0},
  \end{align*}
  the translated $T_0$, so that $\bm\lambda(nP)\in T_\ell$ if and only if $\bm\lambda(nP+\ell)\in T_0$. To ease the notation write $\gamma_i=\lambda_i^P$. Taking subsequences, following \autoref{sec:non-degeneracy}, we divide the $\gamma_i$ into the independent ones and dependent ones. So there exists a partition of $\set{1,\ldots,d}$ into subsets $I,D$, which by rearranging assume that $I:=\set{1,\ldots,\nu}$ and $D:=\set{\nu+1,\ldots, d}$, and rationals $q_{i,j}\in\rat$, $i\in D, j\in I$, such that we can write
  \begin{align*}
    \left(\gamma_1,\ldots,\gamma_d\right)=\left(\gamma_1,\ldots, \gamma_\nu,\prod_{j\in I}\gamma_j^{q_{\nu+1,j}},\ldots, \prod_{j\in I}\gamma_j^{q_{d,j}} \right).
  \end{align*}
  And furthermore there are no multiplicative relations among the $\gamma_1,\ldots, \gamma_\nu$. Imposing these dependencies for the coordinates $\nu+1,\ldots, d$ on the set $T_\ell$, \ie by requiring that
  \begin{align*}
    \gamma_i=\prod_{j\in I}\gamma_j^{q_{i,j}},
  \end{align*}
  for all $i\in D$, we get a new semialgebraic set, denoted $\hat T_\ell$ which is a subset of $\torus^{\nu}$. By definition we have that for all $n\in\nat$, the following equivalences hold
  \begin{align*}
    \left(\gamma_1^n,\ldots,\gamma_\nu^n\right)\in \hat T_\ell\quad \Leftrightarrow \quad (\gamma_1^n,\ldots,\gamma_d^n)\in T_\ell\quad \Leftrightarrow \quad t_\ell(n)=0.
  \end{align*}
  Since $\gamma_1,\ldots,\gamma_\nu$ have no multiplicative relations, applying \autoref{th:cassels}, and the equivalence just above,  it follows that the density of $0$ in the word $t_\ell$ is equal to the Lebesgue measure of $\hat T_\ell$. Since the latter is a semialgebraic set, its Lebesgue measure is nonzero if and only if $\hat T_\ell$ has nonempty interior. Indeed, a nonempty semialgebraic set that has empty interior must be a finite union of hyper-surfaces which have zero volume; for the converse it holds generally that any set with nonempty interior has positive volume.

  Finally to decide whether $\hat T_\ell$ has nonempty interior we write a sentence in the first order logic of reals by saying that there exists some $r>0$ and $r$-ball that is a subset of $\hat T_\ell$ and decide whether it is true by appealing to \autoref{th:complexity of fo rel}. 
\end{proof}
\begin{thm}
    There is a procedure for the following problem. Given as input:
  \begin{itemize}
  \item semialgebraic sets $S_1,\ldots, S_k\subset \torus^d$,
  \item algebraic numbers $\bm\lambda\in\torus^d$,
    
  \item a pattern $w\in 2^{\set{1,\ldots,k}^*}$,
  \item a rational constant $\epsilon>0$
  \end{itemize}
  compute the density of $w$ in the symbolic orbit of $\bm\lambda$ for $S_1,\ldots, S_k$ up to $\epsilon$ additive precision.  
\end{thm}
\begin{proof}
  From the proof of the preceding theorem, it suffices to only estimate the volume of the semialgebraic set $\hat T_\ell$. For this we proceed as in \autoref{sec:comp}. 
\end{proof}

Having shown that we can compute the density of patterns in a toric word, it remains to show that the symbolic orbit of a \lds is similar to some toric word, which we can effectively construct. Indeed we will now prove that the symbolic orbit differs from a toric word in only a subset of indices that have zero density. Intuitively this is because the orbit of the \lds depends primarily on the dominant eigenvalues.

\begin{thm}
  Let $x_0\in\rat^d$, $M\in\rat^{d\times d}$, be a given \lds and $S_1,\ldots, S_k$ semialgebraic subsets of $\rel^d$. Denote by $s$ the symbolic orbit of $(x_0,M)$ for $S_1,\ldots, S_k$. Then there exists a toric word $t$ such that the set of $n\in\nat$ for which
  \begin{align*}
    s(n)\ne t(n),
  \end{align*}
  has density zero. 
\end{thm}
\begin{proof}
  Let $\sq{u_n}$ be a \lrs and denote by $s_1$ the infinite word over the alphabet $\set{p,\tilde p}$, where we put the letter $p$ in position $n$ (\ie $s_1(n)=p$) if and only if $u_n\ge 0$. We claim that:
  \begin{clm}
    \label{cl:claim}
  There exists a toric word $t_1$ over the same alphabet that differs from $s_1$ only in a set of density zero.
\end{clm}
\begin{proof}[Proof of \autoref{cl:claim}]
Let $\alpha_1,\ldots,\alpha_r$ be the dominant characteristic roots of the sequence, \ie those of maximal modulus, assumed distinct. Divide the sequence $\sq{u_n}$ by $|\alpha_1|^nn^{m-1}$ where $m$ is the maximal multiplicity of the characteristic roots that have maximal modulus; same as in~\eqref{eq:vn}, to get a new sequence
  \begin{align*}
    v_n\defeq \underbrace{\sum_{i=1}^rc_i\alpha_i^n}_{D(n)}+R(n), 
  \end{align*}
  where $R(n)$ is some remainder that tends to zero as $n$ grows larger. This new sequence has the exact same signs as the sequence $\sq{u_n}$ and therefore the same $\omega$-word $s_1$. We let $t_1$ be the toric word that is the symbolic orbit of $\bm\alpha$ for $S,\tilde S$, where $S\subset\torus^r$ is the semialgebraic set
  \begin{align*}
    \set{\bm z\in \torus^r\st \sum_{i=1}^rc_iz_i^n\ge 0},
  \end{align*}
  and $\tilde S$ its relative complement in $\torus^r$. So $t_1$ corresponds to the signs of $D(n)$, while $s_1$ corresponds to the signs of $D(n)+R(n)$. By splitting into sub-words $nP+\ell$, $0\le \ell<P$, for $P$ defined in \autoref{sec:non-degeneracy} and applying \autoref{lem:R does not matter}, we get that sub-words differ only on a set of density zero (because the sign of the remainder matters only rarely). Since the union of a finite number of subsets of $\nat$ that have density zero, also has density zero, the claim follows. 
\end{proof}

Let $f$ be a polynomial in $\intg[x_1,\ldots, x_d]$. The sequence
\begin{align*}
  f(x_0M^n), n\in\nat,
\end{align*}
is an \lrs. This is because every component of $x_0M^n,n\in\nat$ is itself an \lrs, and these sequences are closed under point-wise addition and product. It follows that if we denote by $S$ the semialgebraic set
\begin{align*}
  \set{\bm x\in\rel^d\st f(\bm x)\ge 0},
\end{align*}
and by $\tilde S$ its complement, the symbolic orbit of the \lds $(x_0, M)$ for $S,\tilde S$ differs from a toric word only on a set of density zero; due to the claim above, \autoref{cl:claim}.

Suppose that we have two such symbolic orbits $s_1$, $s_2$, ($\omega$-words over the alphabet $\set{p,\tilde p}$ with the semantics above), one for a polynomial $f_1$ and another for another polynomial $f_2$. Let $t_1$ respectively $t_2$ be the toric words to which they are similar. We can take the product $t_1\times t_2$ which is also a toric word (thanks to \autoref{prop:closed products}) and then coarsen it by mapping $(p,p)$, $(p,\tilde p)$, $(\tilde p, p)$ to the same letter, say $a$, and $(\tilde p,\tilde p)$ to the other letter, say $b$. The resulting word is toric (thanks to \autoref{prop:closed under coarsenings}) and it differs in only a set of density zero from the symbolic orbit of $(x_0,M)$ for the union of semialgebraic sets corresponding to $f_1\ge 0$ and $f_2\ge 0$. Similarly we proceed for intersection. Since semialgebraic sets are just unions and intersections of sets of $\bm x$ for which $f(\bm x)\ge 0$, the theorem follows. 
\end{proof}

\bibliographystyle{alphaurl}
\bibliography{bibliography}
\end{document}

%% file: macros.tex
\newcommand{\set}[1]{\left\{#1\right\}}
\newcommand{\st}{\ :\ }
\newcommand{\nat}{\mathbb{N}}
\newcommand{\intg}{\mathbb{Z}}
\newcommand{\rel}{\mathbb{R}}
\newcommand{\rat}{\mathbb{Q}}
\newcommand{\comp}{\mathbb{C}}
\newcommand{\alg}{\overline \rat}

\newcommand{\ii}{\mathbf{i}}

\newcommand{\lrs}{{\sc lrs}\xspace}
\newcommand{\lds}{{\sc lds}\xspace}
\newcommand{\mso}{{\sc mso}\xspace}

\newcommand{\sq}[1]{\langle #1 \rangle_{n\in\nat}}

\newcommand{\ie}{{\em i.e.}\xspace}
\newcommand{\eg}{{\em e.g.}\xspace}

\newcommand\defeq{\mathrel{\overset{\makebox[0pt]{\mbox{\normalfont\tiny\sffamily def}}}{=}}}

\newcommand\dens{\mathcal{D}}

\newcommand{\ptime}{{\sc ptime}\xspace}
\newcommand{\pspace}{{\sc pspace}\xspace}

\newcommand{\np}{{\sc np}\xspace}
\newcommand{\sat}{{\sc sat}\xspace}

\newcommand{\feas}{{\sc feas}\xspace}

\newcommand\torus{\mathbb T}
\newcommand\udens{\hat\dens}
\newcommand\indic{\mathds{1}}
